\documentclass[12pt,reqno,english]{smfart}

\usepackage{amsmath,amssymb,smfthm,stmaryrd,epic,enumerate,dsfont}
\usepackage[english]{babel}
\usepackage[all]{xy}
\usepackage{csquotes}
\usepackage{mathabx}
\usepackage{xcolor}

\usepackage{mathrsfs}
\usepackage[active]{srcltx}

\newcommand*{\htarrow}{\lhook\joinrel\twoheadrightarrow}

\NumberTheoremsIn{subsection}\SwapTheoremNumbers

\begin{document}

\newcommand{\REMARK}[1]{\marginpar{\tiny #1}}
\newtheorem{thm}{Theorem}[section]
\newtheorem{defin}[thm]{Definition}
\newtheorem{Remark}[thm]{Remark}
\numberwithin{equation}{subsection}
\newtheorem{coro}[thm]{Corollary}
\newtheorem{prop}[thm]{Proposition}
\newtheorem{lemma}[thm]{Lemma}

\newtheorem{notas}[thm]{Notations}
\newtheorem{nota}[thm]{Notation}
\newtheorem*{thm*}{Theorem}
\newtheorem*{prop*}{Proposition}
\newtheorem{conje}[thm]{Conjecture}
\newtheorem*{conj*}{Conjecture}
\newtheorem{hyp}[thm]{Hypothesis}

\def\Tm{{\mathbb T}}
\def\Um{{\mathbb U}}
\def\Am{{\mathbb A}}
\def\Fm{{\mathbb F}}
\def\Mm{{\mathbb M}}
\def\Nm{{\mathbb N}}
\def\Pm{{\mathbb P}}
\def\Qm{{\mathbb Q}}
\def\Zm{{\mathbb Z}}
\def\Dm{{\mathbb D}}
\def\Cm{{\mathbb C}}
\def\Rm{{\mathbb R}}
\def\Gm{{\mathbb G}}
\def\Lm{{\mathbb L}}
\def\Km{{\mathbb K}}
\def\Om{{\mathbb O}}
\def\Em{{\mathbb E}}
\def\Xm{{\mathbb X}}

\def\BC{{\mathcal B}}
\def\QC{{\mathcal Q}}
\def\TC{{\mathcal T}}
\def\ZC{{\mathcal Z}}
\def\AC{{\mathcal A}}
\def\CC{{\mathcal C}}
\def\DC{{\mathcal D}}
\def\EC{{\mathcal E}}
\def\FC{{\mathcal F}}
\def\GC{{\mathcal G}}
\def\HC{{\mathcal H}}
\def\IC{{\mathcal I}}
\def\JC{{\mathcal J}}
\def\KC{{\mathcal K}}
\def\LC{{\mathcal L}}
\def\MC{{\mathcal M}}
\def\NC{{\mathcal N}}
\def\OC{{\mathcal O}}
\def\PC{{\mathcal P}}
\def\UC{{\mathcal U}}
\def\VC{{\mathcal V}}
\def\AC{{\mathcal A}}
\def\SC{{\mathcal S}}
\def\RC{{\mathcal R}}

\def\BF{{\mathfrak B}}
\def\AF{{\mathfrak A}}
\def\GF{{\mathfrak G}}
\def\EF{{\mathfrak E}}
\def\CF{{\mathfrak C}}
\def\DF{{\mathfrak D}}
\def\JF{{\mathfrak J}}
\def\LF{{\mathfrak L}}
\def\MF{{\mathfrak M}}
\def\NF{{\mathfrak N}}
\def\XF{{\mathfrak X}}
\def\UF{{\mathfrak U}}
\def\KF{{\mathfrak K}}
\def\FF{{\mathfrak F}}

\def \longmapright#1{\smash{\mathop{\longrightarrow}\limits^{#1}}}
\def \mapright#1{\smash{\mathop{\rightarrow}\limits^{#1}}}
\def \lexp#1#2{\kern \scriptspace \vphantom{#2}^{#1}\kern-\scriptspace#2}
\def \linf#1#2{\kern \scriptspace \vphantom{#2}_{#1}\kern-\scriptspace#2}
\def \linexp#1#2#3 {\kern \scriptspace{#3}_{#1}^{#2} \kern-\scriptspace #3}

\def \Sel {{\mathop{\mathrm{Sel}}\nolimits}}
\def \Ext{\mathop{\mathrm{Ext}}\nolimits}
\def \ad{\mathop{\mathrm{ad}}\nolimits}
\def \sh{\mathop{\mathrm{Sh}}\nolimits}
\def \irr{\mathop{\mathrm{Irr}}\nolimits}
\def \FH{\mathop{\mathrm{FH}}\nolimits}
\def \FPH{\mathop{\mathrm{FPH}}\nolimits}
\def \coh{\mathop{\mathrm{Coh}}\nolimits}
\def \res{\mathop{\mathrm{Res}}\nolimits}
\def \op{\mathop{\mathrm{op}}\nolimits}
\def \rec {\mathop{\mathrm{rec}}\nolimits}
\def \art{\mathop{\mathrm{Art}}\nolimits}
\def \vol {\mathop{\mathrm{vol}}\nolimits}
\def \cusp {\mathop{\mathrm{Cusp}}\nolimits}
\def \scusp {\mathop{\mathrm{Scusp}}\nolimits}
\def \Iw {\mathop{\mathrm{Iw}}\nolimits}
\def \JL {\mathop{\mathrm{JL}}\nolimits}
\def \speh {\mathop{\mathrm{Speh}}\nolimits}
\def \isom {\mathop{\mathrm{Isom}}\nolimits}
\def \Vect {\mathop{\mathrm{Vect}}\nolimits}
\def \groth {\mathop{\mathrm{Groth}}\nolimits}
\def \hom {\mathop{\mathrm{Hom}}\nolimits}
\def \deg {\mathop{\mathrm{deg}}\nolimits}
\def \val {\mathop{\mathrm{val}}\nolimits}
\def \det {\mathop{\mathrm{det}}\nolimits}
\def \rep {\mathop{\mathrm{Rep}}\nolimits}
\def \spec {\mathop{\mathrm{Spec}}\nolimits}
\def \fr {\mathop{\mathrm{Fr}}\nolimits}
\def \frob {\mathop{\mathrm{Frob}}\nolimits}
\def \ker {\mathop{\mathrm{Ker}}\nolimits}
\def \im {\mathop{\mathrm{Im}}\nolimits}
\def \Red {\mathop{\mathrm{Red}}\nolimits}
\def \red {\mathop{\mathrm{red}}\nolimits}
\def \aut {\mathop{\mathrm{Aut}}\nolimits}
\def \diag {\mathop{\mathrm{diag}}\nolimits}
\def \spf {\mathop{\mathrm{Spf}}\nolimits}
\def \Def {\mathop{\mathrm{Def}}\nolimits}
\def \twist {\mathop{\mathrm{Twist}}\nolimits}
\def \supp {\mathop{\mathrm{Supp}}\nolimits}
\def \Id {{\mathop{\mathrm{Id}}\nolimits}}
\def \lie {{\mathop{\mathrm{Lie}}\nolimits}}
\def \Ind{\mathop{\mathrm{Ind}}\nolimits}
\def \ind {\mathop{\mathrm{ind}}\nolimits}
\def \bad {\mathop{\mathrm{Bad}}\nolimits}
\def \top {\mathop{\mathrm{Top}}\nolimits}
\def \ker {\mathop{\mathrm{Ker}}\nolimits}
\def \coker {\mathop{\mathrm{Coker}}\nolimits}
\def \gal {{\mathop{\mathrm{Gal}}\nolimits}}
\def \Nr {{\mathop{\mathrm{Nr}}\nolimits}}
\def \rn {{\mathop{\mathrm{rn}}\nolimits}}
\def \tr {{\mathop{\mathrm{Tr~}}\nolimits}}
\def \Sp {{\mathop{\mathrm{Sp}}\nolimits}}
\def \st {{\mathop{\mathrm{St}}\nolimits}}
\def \sp{{\mathop{\mathrm{Sp}}\nolimits}}
\def \perv{\mathop{\mathrm{Perv}}\nolimits}
\def \tor {{\mathop{\mathrm{Tor}}\nolimits}}
\def \gr {{\mathop{\mathrm{gr}}\nolimits}}
\def \nilp {{\mathop{\mathrm{Nilp}}\nolimits}}
\def \obj {{\mathop{\mathrm{Obj}}\nolimits}}
\def \spl {{\mathop{\mathrm{Spl}}\nolimits}}
\def \unr {{\mathop{\mathrm{Unr}}\nolimits}}
\def \alg {{\mathop{\mathrm{Alg}}\nolimits}}
\def \grr {{\mathop{\mathrm{grr}}\nolimits}}
\def \cogr {{\mathop{\mathrm{cogr}}\nolimits}}
\def \coFil {{\mathop{\mathrm{coFil}}\nolimits}}

\def \rem{{\noindent\textit{Remark.~}}}
\def \rems{{\noindent\textit{Remarks:~}}}
\def \ext {{\mathop{\mathrm{Ext}}\nolimits}}
\def \End {{\mathop{\mathrm{End}}\nolimits}}

\def\semi{\mathrel{>\!\!\!\triangleleft}}
\let \DS=\displaystyle
\def\HT{{\mathop{\mathcal{HT}}\nolimits}}

\def \hi{\HC}
\newcommand*{\tarrow}{\relbar\joinrel\mid\joinrel\twoheadrightarrow}
\newcommand*{\harrow}{\lhook\joinrel\relbar\joinrel\mid\joinrel\rightarrow}
\newcommand*{\rarrow}{\relbar\joinrel\mid\joinrel\rightarrow}
\def \coim {{\mathop{\mathrm{Coim}}\nolimits}}
\def \can {{\mathop{\mathrm{can}}\nolimits}}
\def\LFF{{\mathscr L}}

\setcounter{secnumdepth}{3} \setcounter{tocdepth}{3}

\def \Fil{\mathop{\mathrm{Fil}}\nolimits}
\def \CoFil{\mathop{\mathrm{CoFil}}\nolimits}
\def \Fill{\mathop{\mathrm{Fill}}\nolimits}
\def \CoFill{\mathop{\mathrm{CoFill}}\nolimits}
\def\SF{{\mathfrak S}}
\def\PF{{\mathfrak P}}
\def \EFil{\mathop{\mathrm{EFil}}\nolimits}
\def \EFill{\mathop{\mathrm{EFill}}\nolimits}
\def \FP{\mathop{\mathrm{FP}}\nolimits}

\let \longto=\longrightarrow
\let \oo=\infty

\let \d=\delta
\let \k=\kappa

%\newcounter{num}
\renewcommand{\theequation}{\arabic{section}.\arabic{thm}}
\newcommand{\marque}{\addtocounter{thm}{1}
{\smallskip \noindent \textit{\thethm}~---~}}

\renewcommand\atop[2]{\ensuremath{\genfrac..{0pt}{1}{#1}{#2}}}

\newcommand\atopp[2]{\genfrac{}{}{0pt}{}{#1}{#2}}

\title{Level lowering: a Mazur principle in higher dimension}

%\alttitle{Cohomology of Lubin-Tate spaces: a new natural geometric proof}

\author{Boyer Pascal}
\email{boyer@math.univ-paris13.fr}
\address{Universit\'e Paris 13, Sorbonne Paris Nord \\
LAGA, CNRS, UMR 7539\\ 
F-93430, Villetaneuse (France) \\
Coloss: ANR-19-PRC}

\frontmatter

\begin{abstract}
For a maximal ideal $\mathfrak m$ of some anemic Hecke algebra $\Tm^S_\xi$
of a similitude group of signature $(1,d-1)$, one can
associate a Galois $\overline \Fm_l$-representation 
$\overline \rho_{\mathfrak m}$ 
%which we supposed to be irreducible of dimension $d$, 
as well as a Galois $\Tm_{\xi,\mathfrak m}^S$-representation
$\rho_{\mathfrak m}$. For $l\geq d$, on can also define a monodromy operator
$\overline N_{\mathfrak m}$ as well as 
$N_{\widetilde{\mathfrak m}}$ for every prime ideal
$\widetilde{\mathfrak m} \subset \mathfrak m$,
giving rise to partitions $\underline{\bar d_{\mathfrak m}}$ and
$\underline d_{\widetilde{\mathfrak m}}$ of $d$.
As with Mazur's principle for $GL_2$, analysing the difference between these
partitions, we infer informations about 
%the set of prime ideals $\widetilde{\mathfrak m} \subset \mathfrak m$, i.e. 
the liftings of $\overline \rho_{\mathfrak m}$ in characteristic zero known
as level lowering problem.

%
%We exhibit cases of a level fixing phenomenon for galoisian automorphic  
%representations of a CM field $F$, with dimension $d \geq 2$. 
%The proof rests on the freeness of the localized 
%cohomology groups of KHT Shimura varieties and the strictness of its
%filtration induced by the spectral sequence associated to the filtration of
%stratification of the nearby cycles perverse sheaf at some fixed place $v$ of $F$. 
%The main point is the observation that the action of the unipotent
%monodromy operator at $v$ is then given by those
%on the nearby cycles where its order of nilpotency modulo $l$ equals those in
%characteristic zero. Finally we infer some consequences concerning level raising
%and Ihara's lemma.
\end{abstract}

%% Classification mathématique  (2010)
%\subjclass{11G18, 11G10, 14G35, 14G22, 11F70, 11F80}

%\keywords{Variétés de Shimura, cohomologie de torsion, idéal maximal de l'algèbre de Hecke, 
%localisation de la cohomologie, représentation galoisienne}
%% Mots et expressions clés en anglais :


\maketitle

\pagestyle{headings} \pagenumbering{arabic}

\tableofcontents
%
%\mainmatter
%
%\renewcommand{\theequation}{\arabic{section}.\arabic{subsection}.\arabic{smfthm}}

\section{Introduction}
%\renewcommand{\theequation}{\arabic{equation}}
%\backmatter

Let $F=EF^+$ be a finite CM extension of $\Qm$ with $E/\Qm$ an imaginary
quadratic field and $F^+$ totally real.
Consider then a similitude group $G/\Qm$ as in \S \ref{para-notations}
and a place $v$ of $F$ above a prime number $p=uu^c$ split in $E$ and such
that $G$ is split at $p$ with
$G(\Qm_p)\simeq \Qm_p^\times \times GL_d(F_v) \times \prod_{v \neq w|u} 
(B_w^{op})^\times$, cf. 
\S \ref{para-notations}. For any finite set $S \ni v$ of places of $F$, let $\Tm^S$ 
be the anemic Hecke algebra and, for $\xi$ an algebraic representation of
$G(\Qm)$, we denote by $\Tm^S_\xi$ the quotient of $\Tm^S$ of
$\xi$-cohomological Satake's parameters.
For any  maximal ideal $\mathfrak m$ of $\Tm_\xi^S$, and for a prime ideal
$\widetilde{\mathfrak m} \subset \mathfrak m$,
%of $\Tm^S\otimes_{\overline \Zm_l} \overline \Qm_l$ 
we denote by 
$$\rho_{\widetilde{\mathfrak m}}:\gal_{F,S} \longrightarrow 
GL_d(\overline \Qm_l)$$ 
the Galois $\overline \Qm_l$-representation
associated to $\widetilde{\mathfrak m}$, cf. \cite{h-t}, where
$\gal_{F,S}$ is the Galois group of the maximal extension of $F$ which
is unramified outside $S$.
By Cebotarev's density theorem and the fact that a semi-simple representation
is determined, up to isomorphism, by characteristic polynomials, then
the semi-simple class $\overline{\rho}_{\mathfrak m}$
of the reduction modulo $l$ of 
$\rho_{\widetilde{\mathfrak m}}$ depends only of the maximal
ideal $\mathfrak m$ of $\Tm_\xi^S$ containing $\widetilde{\mathfrak m}$.

\medskip

\noindent \textbf{Main assumptions}: \textit{we now suppose}
\begin{itemize}
\item $[F(\exp(\frac{2i\pi}{l}):F] >d$: \textit{there exists then $v$ as above
such that the order $q_v$ of the residue field at $v$ is of order $>d$ modulo $l$;}

\item \textit{$\overline \rho_{\mathfrak m}$ is absolutely irreducible,}

\item \textit{and 
as a representation of the Weil group at $v$, 
up to the action of the monodromy operator, the semi-simplification of
$\overline \rho_{\mathfrak m,v}$ is a multiplicity free direct sum of characters.}

%
%The set $S_v(\mathfrak m)$ does not contain any subset of the form
%$\{ \lambda,q_v \lambda,\cdots ,q_v^{e_v(l)-1} \lambda \}$ where
%$e_v(l)$ is equal to either the order of $q_v \in \Fm_l^\times$ if it is
%different from $1$ or $l$ otherwise, cf. notation \ref{nota-evl}.}
\end{itemize}

\rem Note first that for every $\widetilde{\mathfrak m} \subset \mathfrak m$,
as $\overline \rho_{\mathfrak m}$ is supposed to be absolutely irreducible, then
$\rho_{\widetilde{\mathfrak m}}$ has, up to homothety, only one stable
$\overline \Zm_l$-lattice.

We are interested in the set $\{ \widetilde{\mathfrak m} \subset \mathfrak m \}$
and the various partitions $\underline{d_{\widetilde{\mathfrak m},v}}$ of $d$
associated to the unipotent operators $N_{\widetilde{\mathfrak m},v}$ 
$$\underline{d_{\widetilde{\mathfrak m},v}}=\Bigl ( n_1(\widetilde{\mathfrak m}) 
\geq n_2(\widetilde{\mathfrak m}) \geq \cdots \geq n_r(\widetilde{\mathfrak m})
\Bigr )$$
with $d=n_1(\widetilde{\mathfrak m})+\cdots n_r(\widetilde{\mathfrak m})$
where the $n_i(\widetilde{\mathfrak m})$ are the sizes of Jordan's blocks of the
nilpotent operator $N_{\widetilde{\mathfrak m}}$. More precisely,
the restriction $\rho_{\widetilde{\mathfrak m},v}$
of $\rho_{\widetilde{\mathfrak m}}$ to the decomposition group at $v$ can be written
as a direct sum
$$\rho_{\widetilde{\mathfrak m},v} \simeq \sp_{n_1(\widetilde{\mathfrak m})}
(\rho_{v,1}) \oplus \cdots \oplus \sp_{n_r(\widetilde{\mathfrak m})}(\rho_{v,r}),$$
where $(n_1(\widetilde{\mathfrak m}) \geq \cdots \geq n_r(\widetilde{\mathfrak m}))$
is a partition of $d$ and
where we suppose the $\rho_{v,i}$ to be characters and
$$\sp_{n_i(\widetilde{\mathfrak m})}(\rho_{v,i})=\rho_{v,i}
(\frac{1-n_i(\widetilde{\mathfrak m})}{2}) \oplus \rho_{v,i}(\frac{3-n_i(\widetilde{\mathfrak m})}{2}) \oplus
\cdots \oplus \rho_{v,i}(\frac{n_i(\widetilde{\mathfrak m})-1}{2}),$$
where $N_{\widetilde{\mathfrak m},v}$ induces isomorphisms 
$\rho_{v,i}(\frac{1-n_i(\widetilde{\mathfrak m})+2\delta}{2}) \longrightarrow
\rho_{v,i}(\frac{1-n_i(\widetilde{\mathfrak m})+2(\delta+1)}{2})$ 
for $0 \leq \delta < n_i(\widetilde{\mathfrak m})-1$ and is trivial on
$\rho_{v,i}(\frac{n_i(\widetilde{\mathfrak m})-1}{2})$.

\begin{nota}
We denote by $T_{\widetilde{\mathfrak m},v}$ the Young diagram of
$$\underline{d_{\widetilde{\mathfrak m},v}}:=(n_1(\widetilde{\mathfrak m}) \geq \cdots \geq n_1(\widetilde{\mathfrak m}))$$ 
labelled by the 
$\rho_{v,i}(\frac{1-n_i(\widetilde{\mathfrak m})+2k}{2})$ so that
$$\Bigl \{ \rho_{v,i}(\frac{1-n_i(\widetilde{\mathfrak m})+2k}{2}): k=0,\cdots,
n_i(\widetilde{\mathfrak m}) -1 \Bigr \}$$ 
are the labels of $i$-th line.
\end{nota}
 
As the order of unipotency of the monodromy operator is trivially less than $d$, 
for $l \geq d$ we can define its 
logarithm in $\overline \Fm_l$ and so define the modulo $l$ nilpotent monodromy
operator $\overline N_{\mathfrak m,v}$ associated to $\overline \rho_{\mathfrak m}$
at the place $v$: recall that as $\overline \rho_{\mathfrak m}$ is supposed to be
irreducible, each of the $\rho_{\widetilde{\mathfrak m}}$ has, up to homothety,
a unique stable $\overline \Zm_l$-lattice so that $\overline N_{\mathfrak m,v}$ is
well defined and does not depend on the choice of $\widetilde{\mathfrak m}
\subset \mathfrak m$. We then denote by $\underline{\bar d_{\mathfrak m,v}}$
the partition of $d$ given by the sizes of Jordan's blocks of
$\overline N_{\mathfrak m,v}$, and $T_{\mathfrak m,v}$ its labelled Young diagram.

\begin{nota}
For any $\widetilde{\mathfrak m} \subset \mathfrak m$, we obtain $T_{\mathfrak m,v}$
from $T_{\widetilde{\mathfrak m},v}$ by breaking into $\delta(\widetilde{\mathfrak m},i)$ 
pieces its $i$-th line. We then denote
by $\delta(\widetilde{\mathfrak m})$ the maximum of the $\delta(\widetilde{\mathfrak m},i)$
and speak about degeneration of monodromy when there exists
$\widetilde{\mathfrak m} \subset \mathfrak m$ with 
$\delta(\widetilde{\mathfrak m}) >0$.
\end{nota}

By classical arguments due to Carayol, we can also define a representation 
$$\rho_{\mathfrak m}: \gal_{F,S} \longrightarrow GL_d(\Tm_{\xi,\mathfrak m}^S),$$
interpolating the $\rho_{\widetilde{\mathfrak m}}$ for all
$\widetilde{\mathfrak m} \subset \mathfrak m$.
By Cebotarev such a $\rho_{\mathfrak m}$ is, up to isomorphism, uniquely 
determined and by construction 
$$\rho_{\mathfrak m} \otimes_{\overline \Zm_l} \overline \Qm_l \simeq
\bigoplus_{\widetilde{\mathfrak m} \subset \mathfrak m} \rho_{\widetilde{\mathfrak m}}$$ 
is semi-simple. 

\medskip

\noindent \textbf{Main results}: \textit{
Under the previous assumptions,}
\begin{itemize}
%\item \textit{the nilpotency order of $N_{\mathfrak m,v} \otimes_{\overline \Zm_l} 
%\overline \Fm_l$ is equal to the maximal of the order of nilpotency of $N_{\widetilde{\mathfrak m},v}$ when $\widetilde{\mathfrak m}$ describes the
%set of prime ideals contained in $\mathfrak m$,}
\item \textbf{Mazur's principle:} \textit{there exists $\widetilde{\mathfrak m} \subset \mathfrak m$ such that
$\underline{d_{\widetilde{\mathfrak m},v}}=\underline{\bar d_{\mathfrak m,v}}$.}

\item \textit{The lenght of $\rho_{\mathfrak m} \otimes_{\overline \Zm_l} \overline \Fm_l$
which is equal to $\sharp \{ \widetilde{\mathfrak m} \subset \mathfrak m \}$,
is greater than $1+\max_{\widetilde{\mathfrak m} \subset \mathfrak m} 
\delta(\widetilde{\mathfrak m}).$}
\end{itemize}
As usual this statement should be translated in its automorphic version. 
\begin{itemize}
\item In the case $d=2$ when there exists a lifting $\widetilde{\mathfrak m}$ with
$\pi_{\widetilde{\mathfrak m},v} \simeq \st_2 (\chi_v)$, our multiplicity free
hypothesis means that $q_v \not \equiv 1 \mod l$. When the modulo
$l$ monodromy operator $\overline N_{\mathfrak m,v}$ is trivial, which means
that the modulo $l$ reduction of $\pi_{\widetilde{\mathfrak m},v}$ is unramified, 
then our statement says that there
exists a lifting $\widetilde{\mathfrak m}'$ of the same level as $\pi_{\widetilde{\mathfrak m}}$ outside $v$ and 
with $\pi_{\widetilde{\mathfrak m}',v}$ unramified. This is exactly the statement of
\textit{the classical Mazur's principle}.

\item For $d \geq 2$ with $\widetilde{\mathfrak m}$ such that 
$\pi_{\widetilde{\mathfrak m},v} \simeq \st_d (\chi_v)$, the multiplicity free hypothesis
is equivalent to asking the order of $q_v \in \Fm_l^\times$ to be greater than $d$.
For $\underline{\bar d_{\mathfrak m,v}}=(t_1 \geq t_2 \geq \cdots \geq t_r)$
the result then insures
the existence of $\widetilde{\mathfrak m}' \subset \mathfrak m$ with
$\pi_{\widetilde{\mathfrak m}' }$ of the same level as $\pi_{\widetilde{\mathfrak m}}$
outside $v$ and such that
$$\pi_{\widetilde{\mathfrak m}',v} \simeq \st_{t_1} (\chi_{v,1}) \times \cdots
\times \st_{t_r}(\chi_{v,r}).$$ 
For any partition $\underline{d}$ of $d$, we denote by $\underline{d}^*$ its
dual partition where the lines of $\underline{d}^*$ are the columns of $\underline{d}$.
Then $\pi_{\widetilde{\mathfrak m}',v}$ has non trivial vectors invariant under
$I_{\underline{d_{\widetilde{\mathfrak m},v}}^*}(\OC_v)$ the parahoric subgroup
associated to the dual partition of $\underline{d_{\widetilde{\mathfrak m},v}}$, cf. 
(\ref{eq-parahoric}).
In the case where 
$$\underline{\bar d_{\mathfrak m,v}}=(1 \geq
\cdots \geq 1),$$ 
$\pi_{\widetilde{\mathfrak m}' ,v}$ is thus unramified. Moreover 
$\{ \widetilde{\mathfrak m} \subset \mathfrak m \}$ is then of order at least $d$.

\item Using the usual automorphic cyclic base change, one can formulate a statement
for automorphic representations $\pi$ of $GL_d$ with $\pi^\vee \simeq \pi^c$.
\end{itemize}

%
%\rem In the end of \S \ref{para-main}, we moreover prove some results on the
%depth of $\rho_{\mathfrak m} \otimes_{\overline \Zm_l} \overline \Fm_l$ in
%the sense of definition \ref{defi-depth} which roughly speaking says that
%if you are able to find $\widetilde{\mathfrak m}_1$ and 
%$\widetilde{\mathfrak m}_2$ such that their monodromy at $v$ are really
%different, then $\rho_{\mathfrak m} \otimes_{\overline \Zm_l} \overline \Fm_l$
%is far from being semi-simple.
%
\medskip

The strategy rests on the study of the cohomology of KHT Shimura vatieties.
More precisely,
for any open compact subgroup $I$ of $G(\Am^\oo)$, 
let denote by $\sh_{I,v} \rightarrow \spec \OC_v$ 
the Kottwitz-Harris-Taylor Shimura variety with level $I$, cf. definition \ref{nota-XI},
where $\OC_v$ is the ring of integers of the local field $F_v$ of $F$ at $v$.
If one believe that Tate conjecture is true, then in the sense of 
definition \ref{defi-typic}, $\mathfrak m$ should be KHT-typic which implies that
$\rho_{\mathfrak m}$
appear in the cohomology of $\sh_{I,v}$ localized at $\mathfrak m$, cf.
proposition \ref{prop-typic}.

\begin{defin} \label{defi-free} (cf. the introduction of \cite{boyer-ihara}) \\
We say that $\mathfrak m$ is KHT-free if the 
cohomology groups of the Kottwitz-Harris-Taylor Shimura variety of notation \ref{nota-XI},
localized at $\mathfrak m$, are free. 
\end{defin}

From \cite{boyer-imj}, any of the following properties
ensure KHT-freeness of $\mathfrak m$.
\begin{itemize}
\item[(1)] There exists a place $w_1 \not \in S$ of $F$ above a prime $p_1$ splits in $E$, 
such that the multi-set
$S_{\mathfrak{m}}(w_1)$ of roots of the characteristic polynomial 
$P_{\mathfrak m,w_1}(X)$ of $\overline \rho_{\mathfrak m}(\frob_{w_1})$,
does not contain any 
sub-multi-set of the shape $\{ \alpha, q_{w_1} \alpha \}$ where $q_{w_1}$ is the order
of the residue field of $F$ at $w_1$. 
This hypothesis is called \emph{generic} in \cite{scholze-cara}.

\item[(2)]  When $[F(\exp(2i\pi/l)):F]>d$, if we suppose the following property to be true,
cf. \cite{boyer-imj} 4.17.
If $\theta: G_F \longrightarrow GL_d(\overline \Qm_l)$ is an irreducible continuous
representation such that for all place $w \not \in S$ above a prime $x \in \Zm$ split in $E$,
then $P_{\mathfrak m,w}(\theta(\frob_w))=0$ 
(resp. $P_{\mathfrak m^\vee,w}(\theta(\frob_w))=0$)
implies that $\theta$ is equivalent to $\overline \rho_{\mathfrak m}$ 
(resp. $\overline \rho_{\mathfrak m^\vee}$), where $\mathfrak m^\vee$ is the maximal
ideal of $\Tm^S_\xi$ associated to the dual mutliset of Satake parameters, cf. \cite{boyer-imj}
notation 4.4. In \cite{emerton-gee}, the authors proved that the previous property 
is verified in each of the following cases:
\begin{itemize}
\item either $\overline \rho_{\mathfrak m}$ is induced from a character of $G_K$ where
$K/F$ is a cyclic galoisian extension;

\item or $l \geq d$ and $SL_d(k) \subset \overline \rho_{\mathfrak m}(G_F) \subset
\overline \Fm_l^\times GL_d(k)$ for some subfield $k \subset \overline \Fm_l$.
\end{itemize}

\item[(3)] In \cite{boyer-jep2}, we proved $\mathfrak m$
to be KHT-free as soon as $\overline \rho_{\mathfrak m}$ is irreducible and
$[F(\exp(2i\pi/l)):F]>d$. 
\end{itemize}

By Cebotarev's theorem, the hypothesis $[F(\exp(2i\pi/l)):F]>d$ allows to
pick places $v$ of $F$ such that the order of $q_v$ modulo $l$ is greater than $d$.
In theorem \ref{theo-main} we then add the assumptions that $\mathfrak m$ is
both KHT-free and
KHT-typic. If one believe in Tate's conjecture then $\mathfrak m$ should always
be KHT-typic, cf. proposition \ref{prop-typic}. 

%
%If $\theta: G_F \longrightarrow GL_d(\overline \Qm_l)$ is an irreducible continuous
%representation sur that for all place $w \not \in S$ above a prime $x \in \Zm$ split in $E$,
%then $P_{\mathfrak m,w}(\theta(\frob_w))=0$ (resp. $P_{\mathfrak m^\vee,w}(\theta(\frob_w))=0$)
%then $\theta$ is equivalent to $\overline \rho_{\mathfrak m}$ 
%(resp. $\overline \rho_{\mathfrak m^\vee}$), where $\mathfrak m^\vee$ is the maximal
%ideal of $\Tm_S$ associated to the dual mutliset of Satake parameters, cf. \cite{boyer-imj}
%notation 4.4.
%
%
%\rem In \cite{emerton-gee}, 

%
%\begin{defin} \label{defi-free}
%We say that $\scusp_v(\mathfrak m)$ is free if for all $\varrho \in \scusp_v(\mathfrak m)$ either
%\begin{itemize}
%\item $\varrho \not \simeq \varrho(1)$ and the intersection of the multiset $\scusp_v(\mathfrak m)$
%with the Zelevinsky line of $\varrho$ is both multiplicity free and not equal to it;
%
%\item when $\varrho \simeq \varrho\{ 1 \}$, we ask $l \geq d+1$.
%\end{itemize}
%\end{defin}
%

The two main ingredients of the proof are:
\begin{itemize}
\item first we show, cf. proposition \ref{prop-torsion}, 
that the filtration given by the monodromy operator at $v$, of
the middle cohomology group of the geometric generic fiber of $\sh_{I,v}$,
is strict, i.e. also gives a filtration of the modulo $l$ cohomology group.

\item We construct a $\overline \Zm_l$-structure of the monodromy operator on
the nearby cycle perverse sheaf which coincides with the usual monodromy operator
over $\overline \Qm_l$ and such that $\overline \Fm_l$ its order of nilpotency is $d$.
\end{itemize}
From this two points, we infer both 
\begin{itemize}
\item a nilpotent monodromy operator 
$\overline N^{coho}_{\mathfrak m,v}$ acting on the $\overline \Fm_l$-cohomology 
group in middle degree of $\sh_{I,v}$;

\item a $\overline \Zm_l$-monodromy operator $N_{\mathfrak m,v}$ on 
$\rho_{\mathfrak m}$,
\end{itemize}
such that the following
\noindent \textbf{main observation} occurs (cf. corollary \ref{coro-main}): \textit{
Suppose that $\mathfrak m$ is KHT-free and KHT-typic in the sense of 
definition \ref{defi-typic}, then the (multi)-set 
of Jordan's blocks of 
$N_{\mathfrak m,v} \otimes_{\overline \Zm_l} \overline \Fm_l$ is the union 
over $\{ \widetilde{\mathfrak m} \subset \mathfrak m \}$, of
the (multi)-set of Jordan's block of $N_{\widetilde{\mathfrak m},v}$.}

\medskip

The existence of various liftings of $\overline \rho_{\mathfrak m}$ with
different levels at $v$ produce constraints on the partition
$\underline{\bar d_{\mathfrak m,v}}$ given by
the Jordan blocks of the modulo $l$ monodromy operator, which have to be smaller 
than all the $\underline{d_{\widetilde{\mathfrak m},v}}$.
 We can then see our result as a reciprocal statement. 
 
\rem Concerning the link with the results of \cite{boyer-mazur}, we refer the reader
to \S \ref{para-reducible} where we consider 
\textbf{the case $\overline \rho_{\mathfrak m}$ reducible.}

\medskip

To see that situations as in the main results really exist, 
we can follow the strategy of Ribet in
\cite{ribet-congruence}, where he considers an absolutely irreducible representation
$$\overline \rho: \gal(\overline \Qm/\Qm) \longrightarrow GL_2(\overline \Fm_l),$$
which is modular of level $N$, meaning that it arises from a cusp form of weight $2$
and trivial character on $\Gamma_0(N)$. Then for a prime $p \nmid lN$
such that
$$\tr \overline \rho(\frob_p) \equiv \pm (p+1) \mod l,$$
with $p+1 \not \equiv 0 \mod l$,
he proves that $\overline \rho$ also arises from a modular form of level $pN$
which is $p$-new, i.e. the automorphic representation associated to this
modular form 
has a local component at $p$ which is isomorphic to the
Steinberg representation of $GL_2(\Qm_p)$. 
%
%b) Using Cebotarev density theorem, there exists 
%infinitely many primes $p \nmid lN$ such that $\overline \rho(\frob_p)$ is 
%conjugate to $\overline \rho(c)$ where $c$ is
%the complex conjugation $c \in \gal(\overline \Qm/\Qm)$. 
%Then we have both 
%$$p \equiv -1 \mod l \hbox{ and } \tr \overline \rho(\frob_p) \equiv 0 \mod l,$$ 
%so that $\overline \rho$ is $p$-new
%for infinitely many $p$. 

In \cite{sor1} Sorensen generalizes this level raising 
congruences in higher dimension for a connected reductive group $G$ over
a totally real field $F^+$ such that $G_\oo$ is compact. One might also
look at \cite{thorne-raising} theorem 1.1 and theorem 4.1, 
for the case of automorphic representations
of unitary type of $GL_{2n}$. More generally it seems that level raising is
more or less settle, cf. theorem 5.1.5, or theorem 4.4.1 of \cite{BLGGT}.

\medskip

We also need to be convinced that such degeneration of the monodromy when
passing modulo $l$, could also appears when $\mathfrak m$ is supposed to
be KHT-free. 
%Here is an example provided  by the anonymous referee whom we thanks to allow us to reproduce it here. 
For example, as pointed out to me, the two following newforms which can be found at
https://www.pnas.org/content/pnas/94/21/11143.full.pdf
\begin{itemize}
\item $f(q)=q-2q^2-q^3+2q^4+q^5+2q^6-2q^7+\cdots \in S_2(\Gamma_0(11))$,

\item $g(q)=q+q^2+2q^3-q^4-2q^5+2q^6-q^7+ \cdots \in S_2(\Gamma_0(77))$,
\end{itemize}
give Galois representations $\rho_f$ and $\rho_g$ which are congruent modulo $3$.
In Serre's 1972 Inventiones paper it is shown that the image of this modulo $3$
representation is $GL_2(\Fm_3)$ so that $\overline \rho_{f}$ is irreducible.
Note then that the $3$-adic representation $\rho_g$ has non trivial
monodromy at $7$ and $\overline \rho_f$ is unramified at $7$ with $7 \not \equiv -1
\mod 3$. Consider then a real quadratic field $F^+$ in which $7$ splits and
which is linearly disjoint over $\Qm$ with the fixed field of $\ker(\overline \rho_g)$.
Let $D/F^+$ a quaternion algebra which is non split t one real place and one
place above $11$. The base change of $g$ to $F^+$ gives a cohomological
automorphic representation of $GL_2(\Am_{F^+})$ whose $3$-adic
Galois representation appears in the cohomology of a Shimura curve attached
to $D$: the $3$-adic monodromy at a place dividing $7$ is then non trivial
while the modulo $3$ monodromy is. Base changing to $F=F^+E$ for
some suitable quadractic imaginary field $E$, we then obtain an example of
a maximal ideal $\mathfrak m$ which is KHT-free and with degeneration
of the monodromy at $v$.

\section{KHT-Shimura varieties and its nearby cycles}
\label{para-notations}

%From now on we fix a maximal ideal $\mathfrak m$ of $\Tm^S$ verifying the hypothesis
%of theorem \ref{theo-principal}. We want to prove that the partition 
%$\underline{d_{\widetilde{\mathfrak m}}}$ associated to any minimal prime ideal
%$\widetilde{\mathfrak m} \subset \mathfrak m$ depends only on $\mathfrak m$.

Let $F=F^+E$ be a CM field with $E/\Qm$ quadratic imaginary and $F^+$ totally real.
Let $B/F$ be a
central division algebra with dimension $d^2$ with an involution of second kind $*$.
For $\beta \in B^{*=-1}$, consider the similitude group $G/\Qm$ defined for any $\Qm$-algebra
$R$ by
$$G(R):=\{ (\lambda,g) \in R^\times \times (B^{op} \otimes_\Qm R)^\times \hbox{ such that }
gg^{\sharp_\beta}=\lambda \},$$
with $B^{op}=B \otimes_{F,c} F$ where $c=*_{|F}$ is the complex conjugation and
$\sharp_\beta$ is the involution $x \mapsto x^{\sharp \beta}:=\beta x^* \beta^{-1}$.
Following \cite{h-t}, we assume from now on that $G(\Rm)$ has signatures $(1,d-1),(0,d),\cdots,(0,d)$.

\begin{defin} \label{defi-spl}
Let $\spl$ be the set of  places $v$ of $F$ such that $p_v:=v_{|\Qm} \neq l$ is split in $E$ and
$B_v^\times \simeq GL_d(F_v)$.  
%For each $I \in \IC$, write $\spl(I)$ for the subset of $\spl$ of places which does not divide the level $I$.
\end{defin}

We now suppose that $p=uu^c$ splits in $E$ so that
$$G(\Qm_p) \simeq \Qm_p^\times \times \prod_{w | u} (B^{op}_w)^\times$$
where $w$ describes the places of $F$ above $u$ and we fix a place
$v \in \spl$ dividing $p$.

\begin{defin}
For a finite set $S$ of places of $\Qm$ containing the places where $G$ is ramified, 
denote by $\Tm^S_{abs}:=\prod_{x \not \in S} \Tm_{x,abs}$ the 
abstract unramified  Hecke algebra
%Hecke algebra where $\unr(I)$ is the union of places $q$ where $G$ is unramified and $I_x$ is 
%maximal, and 
where
$\Tm_{x,abs} \simeq \overline \Zm_l[X^{un}(T_x)]^{W_x}$ for $T_x$ a split torus,
$W_x$ the spherical Weyl group and $X^{un}(T_x)$ is the set of $\overline \Zm_l$-unramified 
characters of $T_x$. 
\end{defin}

\noindent \textit{Example}.
For $x=uu^c$ split in $E$ we have
$$\Tm_{x,abs}=\prod_{w | u} \overline \Zm_l \bigl [T_{w,i}:~ i=1,\cdots,d \bigr ],$$
where $T_{w,i}$ is the characteristic function of
$$GL_d(\OC_w) \diag(\overbrace{\varpi_w,\cdots,\varpi_w}^{i}, \overbrace{1,\cdots,1}^{d-i} ) 
GL_d(\OC_w) \subset  GL_d(F_w).$$

We then denote by $\IC$ the set of open compact subgroups
$$U^p(m_1,\cdots,m_r)=
U^p \times \Zm_p^\times \times \prod_{i=1}^r \ker (\OC_{B_{v_i}}^\times
\longrightarrow (\OC_{B_{v_i}}/\PC_{v_i}^{m_i})^\times )$$
where $U^p$ is any small enough open compact subgroup of $G(\Am^{p,\oo})$ and
$\OC_{B_{v_i}}$ is the maximal order of $B_{v_i}$ with maximal ideal $\PC_{v_i}$ and where
$v=v_1,\cdots,v_r$ are the places of $F$ above $u$ with $p=uu^c$.

\begin{nota}
For $I=U^p(m_1,\cdots,m_r) \in \IC$, we will denote by $I^v(n):=U^p(n,m_2,\cdots,m_r)$.
We also denote by $\spl(I)$ 
the subset of $\spl$ of places which does not divide the level $I$.
\end{nota}

\begin{nota} \label{nota-XI}
As defined in \cite{h-t}, attached to each $I \in \IC$ is a Shimura variety called 
of KHT-type and denoted by 
$$\sh_{I,v} \longrightarrow \spec \OC_v$$ 
where $\OC_v$
denote the ring of integers of the completion $F_v$ of $F$ at $v$. 
\end{nota}

Let $\sigma_0:E \hookrightarrow
\overline{\Qm}_l$ be a fixed embedding and write $\Phi$ for the set of embeddings 
$\sigma:F \hookrightarrow \overline \Qm_l$ whose restriction to $E$ equals $\sigma_0$.
There exists then, cf. \cite{h-t} p.97, an explicit bijection between irreducible algebraic representations 
$\xi$ of $G$ over $\overline \Qm_l$ and $(d+1)$-uple
$\bigl ( a_0, (\overrightarrow{a_\sigma})_{\sigma \in \Phi} \bigr )$
where $a_0 \in \Zm$ and for all $\sigma \in \Phi$, we have $\overrightarrow{a_\sigma}=
(a_{\sigma,1} \leq \cdots \leq a_{\sigma,d} )$. 
We then denote by 
$$V_{\xi,\overline \Zm_l}$$ 
the associated $\overline \Zm_l$-local system on $\sh_{I,v}$.

\begin{nota} \label{nota-TS}
Let $\Tm^S_\xi$ be the image of $\Tm^S_{abs}$ inside 
$$\bigoplus_{i =0}^{2d-2}
\lim_{\atop{\longrightarrow}{I}} H^i_{free}(\sh_{I,\bar \eta_v},V_{\xi,\overline \Zm_l})$$
where the limit concerned the ideals $I$ which are maximal at each places outside $S$, and $\sh_{I,\bar \eta_v}$ is the geometric generic fiber of $\sh_{I,v}$.
\end{nota}

\rem $H^i_{free}$ is the free quotient of the cohomology group $H^i$.
From the main result of \cite{boyer-imj}, the torsion classes of any of the
$H^i(\sh_{I,\bar \eta_v},V_{\xi,\overline \Zm_l})$ raise in characteristic zero, so
one can erase the index \emph{free} in the previous notation.

To each maximal ideal $\widetilde{\mathfrak m}$ of $\Tm^S_\xi[1/l]$,
or equivalently a minimal prime of $\Tm^S_\xi$, which
we now supposed to be non-Eisenstein, is associated an irreducible automorphic 
representation $\Pi_{\widetilde{\mathfrak m}}$
which is $\xi$-cohomological, i.e. there exists
an integer $i$ such that
$$H^i \bigl ( ( \lie ~G(\Rm)) \otimes_\Rm \Cm,U,\Pi_\oo \otimes \xi^\vee \bigr ) \neq (0),$$
where $U$ is a maximal open compact subgroup modulo the center of $G(\Rm)$.
%Let $d_\xi^i(\Pi_\oo)$ be the dimension of this cohomology group.

\begin{nota}
Let denote by $\scusp_v(\widetilde{\mathfrak m})$, the supercuspidal support of
its local component at $v$, denoted $\Pi_{\widetilde{\mathfrak m},v}$. 
Note\footnote{It follows, through the Langlands correspondence, 
from Cebotarev's theorem and the fact that a semi-simple representation is determined, up to isomorphism, by its characteristic
polynomials.}
that the modulo $l$ reduction of $\scusp_v(\widetilde{\mathfrak m})$ is 
independent
of the choice of $\widetilde{\mathfrak m} \subset \mathfrak m$: we denote it
$\scusp_v(\mathfrak m)$. 
\end{nota}

Recall that the geometric special fiber $\sh_{I,\bar s_v}$ of $\sh_{I,v}$, 
is equipped with the Newton stratification 
$$\sh^{\geq d}_{I,\bar s_v} \subset \sh^{\geq d-1}_{I,\bar s_v} \subset \cdots \subset
\sh^{\geq 1}_{I,\bar s_v}=\sh_{I,\bar s_v},$$ 
where for $1 \leq h \leq d$, $\sh_{I,\bar s_v}^{\geq h}$ (resp. $\sh_{I,\bar s_v}^{=h}$)
is the closed (resp. the open) Newton stratum of height $h$ and of pure dimension $d-h$, 
defined as the sub-scheme
where the connected component of the universal Barsotti-Tate group is of rank 
greater or equal to $h$ (resp. equal to $h$).  

Moreover for $1 \leq h < d$, the Newton stratum $\sh_{\IC,\bar s_v}^{=h}$ is 
geometrically induced
under the action of the parabolic subgroup $P_{h,d-h}(F_v)$, defined as the
stabilizer of the first $h$ vectors of the canonical basis of $F_v^d$. Concretely,
cf. \cite{boyer-invent} \S 10.4, this means that there
exists a closed sub-scheme $\sh_{I,\bar s_v,\overline{1_h}}^{=h}$ stabilized by the Hecke 
action of $P_{h,d-h}(F_v)$ and such that
$$\lim_{\atop{\longleftarrow}{n}} \sh_{I^v(n),\bar s_v}^{=h} \simeq\bigl (  
\lim_{\atop{\longleftarrow}{n}}
\sh_{I^v(n),\bar s_v,\overline{1_h}}^{=h}  \bigr )
\times_{P_{h,d-h}(F_v)} GL_d(F_v).$$

\begin{nota}
For a representation $\pi_v$ of $GL_d(F_v)$ with coefficients either 
$\overline \Qm_l$ or $\overline \Fm_l$, and $n \in \frac{1}{2} \Zm$, we set 
$\pi_v \{ n \}:= \pi_v \otimes \nu^n$ where $\nu(g):=q_v^{-\val \det(g)}$. Recall that
the normalized induction of two representations $\pi_{v,1}$ and $\pi_{v,2}$ 
of respectively $GL_{n_1}(F_v)$ and $GL_{n_2}(F_v)$ is
$$\pi_1 \times \pi_2:=\ind_{P_{n_1,n_2}(F_v)}^{GL_{n_1+n_2}(F_v)}
\pi_{v,1} \{ \frac{n_2}{2} \} \otimes \pi_{v,2} \{-\frac{n_1}{2} \},$$
and we define inductively 
$$\pi_1 \times \cdots \times \pi_s:=(\pi_1 \times \cdots \times \pi_{s-1}) \times \pi_s
\simeq \pi_1 \times (\pi_2 \times \cdots \times \pi_s).$$
\end{nota}

Recall that a representation
$\pi_v$ of $GL_d(F_v)$ is called \emph{cuspidal} (resp. \emph{supercuspidal})
if it is not a subspace (resp. a subquotient) of a proper parabolic induced representation.
When the field of coefficients is of characteristic zero then these two notions coincides,
but this is no more true over $\overline \Fm_l$.
For example the modulo $l$ reduction of an irreducible 
$\overline \Qm_l$-representation is still cuspidal but not necessary supercuspidal,
its supercuspidal support being a Zelevinsky segment.
%$[\varrho,\varrho \nu^{s-1}]$ where $\varrho$ is irreducible supercuspidal and $s$ is either equal to $1$ or of the following shape $s=m(\varrho)l^k$ for $k \in \Nm$. 

% 
%\begin{defin} \label{defi-type}
%In the former case we say that $\pi_v$ is of $\varrho$-type $-1$ and otherwise of 
%$\varrho$-type $k$.
%\end{defin}

\begin{defin} \label{defi-LT}
\label{defi-rep} (see \cite{zelevinski2} \S 9 and \cite{boyer-compositio} \S 1.4)
Let $g$ be a divisor of $d=sg$ and $\pi_v$ an irreducible cuspidal 
$\overline \Qm_l$-representation of $GL_g(F_v)$. 
The induced representation
\begin{equation} \label{eq-ind-rep}
\pi_v\{ \frac{1-s}{2} \} \times \pi_v \{ \frac{3-s}{2} \} \times \cdots \times \pi_v \{ \frac{s-1}{2} \}
\end{equation}
holds a unique irreducible quotient (resp. subspace) denoted $\st_s(\pi_v)$ (resp.
$\speh_s(\pi_v)$); it is a generalized Steinberg (resp. Speh) representation.
\end{defin}

\rem For $\chi_v$ a character, $\speh_s(\chi_v)$ is the character $\chi_v \circ \det$ of $GL_s(F_v)$.

Let $\pi_v$ be an irreducible cuspidal $\overline \Qm_l$-representation of $GL_g(F_v)$ and
fix $t \geq 1$ such that $tg \leq d$. Thanks to Igusa varieties, Harris and Taylor
constructed a local system on  $\sh^{=tg}_{\IC,\bar s_v,\overline{1_{tg}}}$
$$\LC_{\overline \Qm_l}(\pi_v[t]_D)_{\overline{1_{tg}}}=\bigoplus_{i=1}^{e_{\pi_v}} 
\LC_{\overline \Qm_l}(\rho_{v,i})_{\overline{1_{tg}}}$$
where 
\begin{itemize}
\item $\pi_v[t]_D$ is the representation of $D_{v,tg}^\times$ which
is the image of the contragredient of $\st_t(\pi_v)$ by the Jacquet-Langlands correspondence,

\item $D_{v,tg}$ is the central division algebra over $F_v$ with
invariant $tg$,

\item with maximal order denoted by $\DC_{v,tg}$ 

\item and with
$(\pi_v[t]_D)_{|\DC_{v,tg}^\times}=\bigoplus_{i=1}^{e_{\pi_v}} \rho_{v,i}$ 
with $\rho_{v,i}$ irreductible. 
\end{itemize}
The Hecke action of $P_{tg,d-tg}(F_v)$ is then given through its quotient 
$$P_{tg,d-tg}(F_v) \twoheadrightarrow GL_{tg}(F_v) \times GL_{d-tg}(F_v)
\twoheadrightarrow GL_{d-tg}(F_v) \times \Zm,$$ 
where $GL_{tg}(F_v) \times GL_{d-tg}(F_v)$ is the Levi quotient of the parabolic
$P_{tg,d-tg}(F_v)$ and the second map is given by the valuation of the determinant map
$GL_{tg}(F_v) \twoheadrightarrow \Zm$.
These local systems have stable
$\overline \Zm_l$-lattices and we will write simply $\LC(\pi_v[t]_D)_{\overline{1_{tg}}}$
for any $\overline \Zm_l$-stable lattice that we do not want to specify.

\begin{notas} For $\Pi_t$ any $\overline \Qm_l$-representation of 
$GL_{tg}(F_v)$, and 
$\Xi:\frac{1}{2} \Zm \longrightarrow \overline \Zm_l^\times$ defined by 
$\Xi(\frac{1}{2})=q^{1/2}$, we introduce
$$\widetilde{HT}_{\overline \Qm_l,\overline{1_{tg}}}(\pi_v,\Pi_t):=
\LC_{\overline \Qm_l}(\pi_v[t]_D)_{\overline{1_{tg}}} 
\otimes \Pi_t \otimes \Xi^{\frac{tg-d}{2}}$$
and its induced version
$$\widetilde{HT}_{\overline \Qm_l}(\pi_v,\Pi_t):=\Bigl ( \LC_{\overline \Qm_l}(\pi_v[t]_D)_{\overline{1_{tg}}} 
\otimes \Pi_t \otimes \Xi^{\frac{tg-d}{2}} \Bigr) \times_{P_{tg,d-tg}(F_v)} GL_d(F_v),$$
where the unipotent radical of $P_{tg,d-tg}(F_v)$ acts trivially and the action of
$$(g^{\oo,v},\left ( \begin{array}{cc} g_v^c & * \\ 0 & g_v^{et} \end{array} \right ),\sigma_v) 
\in G(\Am^{\oo,v}) \times P_{tg,d-tg}(F_v) \times W_v$$ 
where $W_v$ is the Weil group at $v$, is given
\begin{itemize}
\item by the action of $g_v^c$ on $\Pi_t$ and 
$\deg(\sigma_v) \in \Zm$ on $ \Xi^{\frac{tg-d}{2}}$, where
$\deg: W_v \longrightarrow \Zm$ sends geometric Frobenius to $1$,

\item and the action of $(g^{\oo,v},g_v^{et},\val(\det g_v^c)-\deg \sigma_v)
\in G(\Am^{\oo,v}) \times GL_{d-tg}(F_v) \times \Zm$ on $\LC_{\overline \Qm_l}
(\pi_v[t]_D)_{\overline{1_{tg}}} \otimes \Xi^{\frac{tg-d}{2}}$.
\end{itemize}
We also introduce
$$HT_{\overline \Qm_l,\overline{1_{tg}}}(\pi_v,\Pi_t):=
\widetilde{HT}_{\overline \Qm_l,\overline{1_{tg}}}(\pi_v,\Pi_t)[d-tg],$$
and the perverse sheaf
$$P_{\overline \Qm_l}(t,\pi_v)_{\overline{1_{tg}}}:=j^{=tg}_{\overline{1_{tg}},!*} HT_{\overline \Qm_l,\overline{1_{tg}}}(\pi_v,\st_t(\pi_v))
\otimes \Lm(\pi_v),$$
and their induced version, $HT_{\overline \Qm_l}(\pi_v,\Pi_t)$ and 
$P_{\overline \Qm_l}(t,\pi_v)$, where 
$$j^{=tg}=i^{tg} \circ j^{\geq tg}:\sh^{=tg}_{\IC,\bar s_v} \hookrightarrow
\sh^{\geq tg}_{\IC,\bar s_v} \hookrightarrow \sh_{\IC,\bar s_v}$$ 
and $\Lm$ is the local Langlands correspondence composed by contragredient.

We will also denote by $HT_{\overline \Qm_l,\xi}(\pi_v,\Pi_t):= HT_{\overline \Qm_l}(\pi_v,\Pi_t) \otimes V_\xi$ 
and similarly for the other notations as for example
$P_{\overline \Qm_l,\xi}(t,\pi_v):=P_{\overline \Qm_l}(t,\pi_v) \otimes V_\xi$.

Finally we denote by $\sh^{=h}_{\IC,\bar s_v,\neq 1}:=\sh^{=h}_{\IC,\bar s_v}
\setminus \sh^{=h}_{\IC,\bar s_v,\overline 1_h}$ and
$$j^{=h}_{\neq 1}:\sh^{=h}_{\IC,\bar s_v,\neq 1} \hookrightarrow
\sh^{=h}_{\IC,\bar s_v} \hookrightarrow \sh_{\IC,\bar s_v}.$$
\end{notas}

\rems 
\begin{itemize}
\item We will simply denote by $P(t,\pi_v)$ any
$\overline \Zm_l$-lattice of $P_{\overline \Qm_l}(t,\pi_v)$ that we do not
want to precise except that it is stable under the various actions. We will
use a similar convention for the other sheaves introduced before.
When considering $\overline \Fm_l$-coefficients, we will put 
$\overline \Fm_l$ in place of $\overline \Qm_l$ in the notations.

\item Recall that $\pi'_v$ is said inertially equivalent to $\pi_v$, and we write
$\pi_v \sim_i \pi'_v$,
if there exists a character $\zeta: \Zm \longrightarrow  \overline \Qm_l^\times$ such that
$\pi'_v \simeq \pi_v \otimes (\zeta \circ \val \circ \det)$. We denote by
$e_{\pi_v}$ the order of the inertial class of $\pi_v$.

\item Note, cf. \cite{boyer-invent2} 2.1.4, that $P_{\overline \Qm_l}(t,\pi_v)$ depends 
only on the inertial class of $\pi_v$ and
$$P_{\overline \Qm_l}(t,\pi_v)=e_{\pi_v} \PC_{\overline \Qm_l}(t,\pi_v)$$ 
where $\PC_{\overline \Qm_l}(t,\pi_v)$ is an irreducible perverse sheaf.

\item Over $\overline \Zm_l$, we also have the $p+$-perverse structure which is dual 
to the usual $p$-structure.
\end{itemize}

\begin{nota} \label{nota-evl}
Let denote by
$$e_v(l)=\left \{ \begin{array}{ll} 
l, & \hbox{if } q_v \equiv 1 \mod l \\
\hbox{the order of $q_v$ modulo } l, & \hbox{otherwise}
\end{array} \right.$$
\end{nota}

\rem For a character\footnote{For a general supercuspidal representation $\pi_v$
whose modulo $l$ reduction $\varrho$ is still supercuspidal, the same is true
if $t < m(\varrho)$ where $m(\varrho)$ is either the order of the Zelevinsky line of
$\varrho$ if it is $>1$, otherwise $m(\varrho):=l$.}
$\chi_v$ and when $t < e_v(l)$, up to homothety there is only one stable 
$\overline \Zm_l$-stable
lattice of $\LC(\pi_v[t]_D)$. From the description of the modulo $l$ reduction of 
$\st_t(\chi_v)$ in \cite{boyer-repmodl}, the same is then
true for $\PC(t,\chi_v)$.

\begin{nota}
Let denote by
$$\Psi_v:=R\Psi_{\eta_v}(\overline \Zm_l[d-1])(\frac{d-1}{2})$$
the nearby cycles autodual free perverse sheaf on the geometric special fiber $\sh_{I,\bar s_v}$
of $\sh_{I,v}$. 
\end{nota}
In cite \cite{boyer-duke} proposition 3.1.3, we proved the following splitting 
$$\Psi_v \simeq \bigoplus_{g=1}^d \bigoplus_{\varrho \in \scusp_{\overline \Fm_l}(g)}
\Psi_\varrho,$$
where $\scusp_{\overline \Fm_l}(g)$ is the set of inertial equivalence classes of
irreducible $\overline \Fm_l$-supercuspidal of $GL_g(F_v)$ and where
$$\Psi_\varrho \otimes_{\overline \Zm_l} \overline \Qm_l \simeq 
\bigoplus_{\pi_v \in \cusp(\varrho)} \Psi_{\pi_v},$$
where $\cusp(\varrho)$ is the set of inertial equivalence classes of irreducible
$\overline \Qm_l$-cuspidal representations which modulo $l$ reduction is inertially
equivalent to $\varrho$.

\noindent \textit{We now fix a $\overline \Fm_l$-character $\varrho$.}
Following the constructions of \cite{boyer-torsion} \S 2.3, 
\begin{itemize}
\item using the adjonction $j_!j^* \rightarrow \Id$, we can  first define a filtration 
$$\Fil^1_!(\Psi_\varrho) \hookrightarrow \cdots \hookrightarrow
\Fil^d_!(\Psi_\varrho)=\Psi_\varrho,$$
where $\Fil^h_!(\Psi_\varrho)$ is the saturated image of
$j^{=h}_! j^{=h,*} \Psi_\varrho \longrightarrow \Psi_\varrho$.
The graded parts $\gr^k_!(\Psi_\varrho)$ are free perverse sheaves. 
Over $\overline \Qm_l$, this filtration coincides with the iterated kernel of $N_v$,
i.e. $\Fil^k_!(\Psi_\varrho) \otimes_{\overline \Zm_l} \overline \Fm_l \simeq
\ker (N_v^k \otimes_{\overline \Zm_l} \overline \Fm_l)$. We also write
$\coFil_!^k(\Psi_\varrho):=\Psi_\varrho/\Fil^k_!(\Psi_\varrho)$.

\item For every $1 \leq h <d$, we have a short exact sequence of free Hecke-perverse
sheaves
\begin{equation} \label{eq-gr!}
0 \rightarrow j^{=h}_{\neq 1,!} j^{=h,*}_{\neq 1} \gr^h_!(\Psi_\varrho) \longrightarrow
\gr^h_!(\Psi_\varrho) \longrightarrow \lexp p j^{=h}_{\overline{1_h},!*} j^{=h,*}_{\overline{1_h}} \gr^h_!(\Psi_\varrho) \rightarrow 0,
\end{equation}
where 
$$j^{=h,*} \gr^h_!(\Psi_\varrho) \otimes_{\overline \Zm_l} \overline \Qm_l
\simeq \bigoplus_{\chi_v \in \cusp(\varrho)} HT_{\overline \Qm_l}(\chi_v,\st_h(\chi_v))
(\frac{1-h}{2}).$$

\item Dually using the adjunction $\Id \rightarrow j_*j^*$, we can define a cofiltration 
$$\coFil^d_*(\Psi_\varrho) \twoheadrightarrow \cdots \twoheadrightarrow 
\coFil^1_*(\Psi_\varrho),$$
where $\coFil^h_*(\Psi_\varrho)$ is the saturated coimage of $\Psi_\varrho
\longrightarrow j^{=h}_* j^{=h,*} \Psi_\varrho$. The graded parts $\gr^h_*(\Psi_\varrho)$
are also given by
$$0 \rightarrow \lexp p j^{=h}_{1,!*} j^{=h,*}_1 \gr^h_*(\Psi_\varrho) \longrightarrow
\gr^h_*(\Psi_\varrho) \longrightarrow
\lexp p j^{=h}_{\neq 1,!*} j^{=h,*}_{\neq 1} \gr^h_*(\Psi_\varrho) \rightarrow 0,$$
where
$$j^{=h,*}_{\neq 1} \gr^h_*(\Psi_\varrho) \otimes_{\overline \Zm_l} \overline \Qm_l
\simeq \bigoplus_{\chi_v \in \cusp(\varrho)} HT_{\overline \Qm_l}(\chi_v,\st_h(\chi_v))(\frac{h-1}{2}).$$

\item We can also refine the previous filtration to obtain $\Fill^\bullet(\Psi_\varrho)$
 whose graded parts
$\grr^r(\Psi_v)$ are free $\overline \Zm_l$-perverse sheaves, cf. \cite{boyer-torsion} \S 1, 
such that
$$\grr^r(\Psi_\varrho) \otimes_{\overline \Zm_l} \overline \Qm_l \simeq 
\bigoplus_{\chi_v \in \cusp(\varrho)}
\lexp p j^{=tg}_{!*} HT (\chi_v,\st_t(\chi_v))(\frac{1-t+2\delta}{2})$$
for some $0 \leq \delta \leq t-1$.
\end{itemize}

\rem As the order of $q_v$ modulo $l$ is $>d$ then $\cusp(\varrho)$ contains only
characters so that, cf. (\ref{eq-ext0}), the Harris-Taylor local systems have only one
intermediate extension, i.e. 
$$\lexp p j^{=h}_{\overline{1_h},!*} j^{=h,*}_{\overline{1_h}} \gr^h_!(\Psi_\varrho) \simeq
\lexp {p+} j^{=h}_{\overline{1_h},!*} j^{=h,*}_{\overline{1_h}} \gr^h_!(\Psi_\varrho).$$

\bigskip

\noindent \textit{Exchange basic step}:
to go from filtration to another, one can repeat the following process to exchange
the order of appearance of two consecutive subquotient:
$$\xymatrix{
& P_1' \ar@{^{(}->}[d] \ar@{^{(}->}[dr] \\
P_2  \ar@{^{(}->}[dr] \ar@{^{(}->}[r] & X \ar@{->>}[d] \ar@{->>}[r] & P_1 \ar@{->>}[dr]  \\
& P_2' \ar@{->>}[dr] & & T \\
& & T, \ar@{=}[ur]
}$$
where 
\begin{itemize}
\item $P_1$ and $P_2$ are two consecutive subquotient in a given filtration
and $X$ is the subquotient gathering them as a subquotient of this filtration.

\item Over $\overline \Qm_l$, the extension $X \otimes_{\overline \Zm_l} \overline \Qm_l$
is split, so that on can write $X$ as an extension of $P'_2$ by $P'_1$ with
$P'_1 \hookrightarrow P_1$ and $P_2 \hookrightarrow P'_2$ have the same cokernel
$T$, a perverse sheaf of torsion.
\end{itemize}

\rem In the particular case when $P_1$ and $P_2$ are intermediate extensions 
of local systems 
living on different strata such that the two associated intermediate extensions
for the $p$ and $p+$ $t$-structure are isomorphic, then $T$ is necessary zero and
$X$ is then split over $\overline \Zm_l$.

Repeating exchange basic steps, on can then pass from $\Fil^\bullet_*(\Psi_\varrho)$
to $\Fil^\bullet_!(\Psi_\varrho)$.

\begin{lemma} \label{lem-socle}
The socle (resp. the cosocle) of 
$\Psi_\varrho \otimes_{\overline \Zm_l} \overline \Fm_l$ is, up to multiplicities,
$j^{=d}_{!*} HT(\varrho,\st_d(\varrho))(\frac{d-1}{2})$
(resp. $j^{=d}_{!*} HT(\varrho,\st_d(\varrho))(\frac{1-d}{2})$).
\end{lemma}

\begin{proof}
The result follows quite immediately from \cite{boyer-duke} where we described
the sheaves of cohomology of the $\Psi_\varrho \otimes_{\overline \Zm_l}
\overline \Fm_l$
as the modulo $l$ reduction of its $\overline \Zm_l$-cohomology sheaves.
From this computation we deduce that the socle of 
 $\gr^h_!(\Psi_\varrho) \otimes_{\overline \Zm_l} \overline \Qm_l$, up to multiplicities,
 is $j^{=d}_{!*} HT_{\overline \Fm_l}(\varrho,\st_d(\varrho))(\frac{d+1-2h}{2})$.
  From the filtration $\Fil^\bullet_!(\Psi_\varrho)$, we then deduce that the socle
 of $\Psi_\varrho \otimes_{\overline \Zm_l} \overline \Fm_l$ 
 \begin{itemize}
 \item contains
 $j^{=d}_{!*} HT(\varrho,\st_d(\varrho))(\frac{d-1}{2})$,
 
 \item and if it contains something else, without considering multiplicities, it has to be 
 $j^{=d}_{!*} HT(\varrho,\st_d(\varrho))(\frac{d-1-2\delta}{2})$ for $1 \leq \delta < d$.
\end{itemize}
We argue by contradiction by assuming that  
$j^{=d}_{!*} HT(\varrho,\st_d(\varrho))(\frac{d-1-2\delta}{2})$ for some $1 \leq \delta < d$,
belongs to the socle. By duality then
$j^{=d}_{!*} HT(\varrho,\st_d(\varrho))(\frac{1-d+2\delta}{2})$ belongs to the cosocle
so that  $HT(\varrho,\st_d(\varrho))(\frac{d-1-2\delta}{2})$ is a subquotient
of $h^0 \Psi_\varrho \otimes_{\overline \Zm_l} \overline \Fm_l$ which contradicts
the main result of \cite{boyer-duke}.

For the cosocle, we conclude by duality.
\end{proof}

\rem The same arguments, with more precautions, 
allow to show that non split extensions inside $\Psi_\varrho
\otimes_{\overline \Zm_l} \overline \Qm_l$ between Harris-Taylor perverse sheaves
remains non split in $\Psi_\varrho \otimes_{\overline \Zm_l}
\overline \Fm_l$. 

%
%\item The above filtration is compatible with the nilpotent
%monodromy operator $N_v$, i.e. for any $r$ the image of 
%$\Fill^r(\Psi_v) \otimes_{\overline \Zm_l} \overline \Qm_l$ under $N_v$ 
%is some $\Fill^{\phi(r)}(\Psi_v)
%\otimes_{\overline \Zm_l} \overline \Qm_l$ for some decreasing function $\phi$.
%

\bigskip

When dealing with sheaves, there is no need to introduce the local system
$V_{\xi,\overline \Zm_l}$ because it suffices to add $\otimes_{\overline \Zm_l}
V_{\xi,\overline \Zm_l}$ to the formulas.
We now consider a fixed local system $V_{\xi,\overline \Zm_l}$ and, following
previous notations, we write $\Psi_{\varrho,\xi}:=\Psi_\varrho \otimes_{\overline \Zm_l}
V_{\xi,\overline \Zm_l}$.
We then have a spectral sequence
\addtocounter{thm}{1}
\begin{equation} \label{eq-ss-gr}
E_1^{p,q}=H^{p+q}(\sh_{I,\bar s_v},\grr^{-p}(\Psi_{v,\xi})) \Rightarrow 
H^{p+q}(\sh_{I,\bar \eta_v},V_{\xi,\overline \Zm_l}).
\end{equation}

As pointed out in \cite{boyer-mazur}, if for some $\mathfrak m$
the spectral sequence is concentrated in middle degree, i.e. $E_{1,\mathfrak m}^{p,q}=0$
for $p+q \neq d-1$, and all the $E_{1,\mathfrak m}^{p,d-1-p}$ are free, then,
for $l > d$,  the action
of the monodromy operator $N_{v,\mathfrak m}^{coho}$ on 
$H^{d-1}(\sh_{I,\bar \eta_v},V_{\xi,\overline \Qm_l})_{\mathfrak m}$
comes from the action of $N_v$ on $\Psi_v \otimes_{\overline \Zm_l} \overline \Qm_l$.
The aim of the next section is to prove that this property 
remains true over $\overline \Fm_l$.

\section{A saturated filtration of the cohomology}

The aim of this section is the following proposition.

\begin{prop} \label{prop-torsion}
Consider a maximal ideal $\mathfrak m$ of $\Tm^S_\xi$ such that:
\begin{itemize}
\item $\overline \rho_{\mathfrak m}$ is irreducible;

\item $\mathfrak m$ is KHT-free;

\item the restriction $\rho_{\mathfrak m,v}$ to the decomposition group at $v$ 
is, in the Grothendieck group and up to the action of $N_v$, the direct sum of character.
We moreover suppose that the set $S_v(\mathfrak m)$ of modulo $l$ eigenvalues of 
$\overline \rho_{\mathfrak m}(\frob_v)$ does not contain any subset of the form
$\{ \lambda,q_v \lambda,\cdots,q_v^{e_v(l)-1} \lambda \}$, where $e_v(l)$ is defined
in \ref{nota-evl}.
\end{itemize}
Then the $E_{1,\mathfrak m}^{p,q}$ 
are torsion free and trivial for $p+q \neq d-1$.
\end{prop}

Note that, as (\ref{eq-ss-gr}) after localization at $\mathfrak m$, 
degenerates at $E_1$ over $\Qm_l$, then the spectral sequence gives us
a saturated filtration of 
$H^{d-1}(\sh_{I,\bar \eta_v},V_{\xi,\overline \Zm_l})_{\mathfrak m}$.
The proof uses Grothendieck-Verdier duality which explains we need that the $p$ and $p+$ intermediate extensions of Harris-Taylor local
systems $HT(\pi_v,\st_t(\pi_v))$ to be isomorphic. 
In \cite{boyer-duke}, for a irreducible cuspidal representation $\pi_v$ of
$GL_g(F_v)$ and $1 \leq t \leq d/g$ such that the modulo $l$ reduction of $\st_t(\pi_v)$
remains irreducible, we proved that
\addtocounter{thm}{1}
\begin{equation} \label{eq-ext0}
\lexp p j^{=tg}_{!*}  HT(\pi_v,\Pi_h) \simeq \lexp {p+} j^{=tg}_{!*} HT(\pi_v,\Pi_h).
\end{equation}
The proof is rather difficult but almost obvious in the case where $\pi_v$
is a character in which case the condition is that $t < e_v(l)$ which is clearly true
with the hypothesis $e_v(l)>d$.

\begin{proof} %(of proposition \ref{prop-torsion}) \\
As $\mathfrak m$ is supposed to be KHT-free, then all the $E_{\oo,\mathfrak m}^n$ are free.
Moreover, as $\overline \rho_{\mathfrak m}$ is irreducible, then,
cf. \cite{boyer-compositio} \S 3.6, the
$E_{1,\mathfrak m}^{p,q} \otimes_{\overline \Zm_l} \overline \Qm_l$ 
%and $E^{p,q}_{\oo,\mathfrak m} \otimes_{\overline \Zm_l} \overline \Qm_l$
are all zero if $p+q \neq d-1$.
As by hypothesis $\rho_{\mathfrak m,v}$ is made of characters, we can consider
direct factors $\Psi_\varrho$ for $\varrho \in \scusp_v(\mathfrak m)$  a character.
Moreover as $e_v(l)>d$, in
$\Psi_\varrho \otimes_{\overline \Zm_l} \overline \Qm_l$ we have only
to deal with characters $\chi_v$ so that,
%with irreducible $\overline \Qm_l$-cuspidal representations $\pi_{v,-1}$ of $\varrho$-type $-1$, or more precisely, if we do not want to use the results of \cite{boyer-duke},with characters $\chi_v$. 
by (\ref{eq-ext0}), the $p$ and $p+$
intermediate extensions coincide.

\begin{prop} \label{prop-resolution1} (cf. \cite{boyer-duke} \S 2.3)
We have the following equivariant resolution
\addtocounter{thm}{1}
\begin{multline} \label{eq-resolution0}
0 \rightarrow j_!^{=d} HT(\chi_v,\st_h(\chi_v \{ \frac{h-d}{2} \} ) \times 
\speh_{d-h}(\chi_v\{ h/2 \} ))
 \otimes \Xi^{\frac{d-h}{2}} \longrightarrow \cdots  \\
\longrightarrow j_!^{=h+1} HT(\chi_v,\st_h(\chi_v (-1/2)) \times \chi_v \{ h/2 \} ) 
\otimes \Xi^{\frac{1}{2}} \longrightarrow \\ j_!^{=h} HT(\chi_v,\st_h(\chi_v)) 
\longrightarrow  \lexp p j_{!*}^{=h} HT(\chi_v,\st_h(\chi_v)) \rightarrow 0.
\end{multline}
\end{prop}

Note that
\begin{itemize}
\item as this resolution is equivalent to the computation of the sheaves 
cohomology groups
of $\lexp p j_{!*}^{=h} HT(\chi_v,\st_h(\chi_v)) $ as explained for example in
\cite{boyer-duke} proposition B.1.5 of appendice B, then, 
over $\overline \Qm_l$, it follows from the main results of \cite{boyer-invent2}.

\item Over $\overline \Zm_l$, as every terms are free perverse sheaves, then
all the maps are necessary strict. 

\item This resolution, for a a general supercuspidal representation
with supercuspidal modulo $l$ reduction, is one of the main result 
of \cite{boyer-duke} \S 2.3. However the case  of a character $\chi_v$ as above,
is almost obvious. Indeed as the strata
$\sh^{\geq h}_{I^v,\bar s_v,1}$ are smooth, then the constant sheaf, up to shift, is perverse and
so equals to the intermediate extension of the constant sheaf, shifted by $d-h$, 
on $\sh^{=h}_{I^v,\bar s_v,1}$. In particular its sheaves cohomology groups
are well known so that the resolution is completely obvious for 
$\lexp p j_{\overline{1_h},!*}^{=h} HT_{\overline{1_h}}(\chi_v,\st_h(\chi_v))$
if one remember that $\speh_i(\chi_v)$ is just the character $\chi_v \circ \det$
of $GL_i(F_v)$.

The stated resolution is then simply the induced version of the resolution
of $\lexp p j^{=h}_{\overline{1_h},!*} HT_{\overline{1_h}}(\chi_v,\st_h(\chi_v))$:
recall that a direct sum of intermediate extensions is still an intermediate 
extension.
\end{itemize}
%
%\begin{proof}
%For the case of a character $\chi_v$ as above, the argument is almost 
%obvious. Indeed as the strata
%$\sh^{\geq h}_{I^v,\bar s_v,1}$ are smooth, then, cf. the proof of the lemma
%\ref{lem-ext0}, the constant sheaf, up to shift, is perverse and
%so equals to the intermediate extension of the constant sheaf, shifted by $d-h$, 
%on $\sh^{=h}_{I^v,\bar s_v,1}$. In particular its sheaves cohomology groups
%are well known so that the resolution is completely obvious for 
%$\lexp p j_{\overline{1_h},!*}^{=h} HT_{\overline{1_h}}(\chi_v,\st_h(\chi_v))$
%if one remember that $\speh_i(\chi_v)$ is just the character $\chi_v \circ \det$
%of $GL_i(F_v)$.
%
%The stated resolution is then simply the induced version of the resolution
%of $\lexp p j^{=h}_{\overline{1_h},!*} HT_{\overline{1_h}}(\chi_v,\st_h(\chi_v))$:
%recall that a direct sum of intermediate extensions is still an intermediate 
%extension.
%% 
%is also the intermediate extension of $HT_{\overline{1_h}}(\chi_v,\st_h(\chi_v))$
%to the whole stratum, and so
%$\Bigl ( \lexp p j^{=h}_{\overline{1_h},!*} HT_{\overline{1_h}}(\chi_v,\st_h(\chi_v))  
%\Bigr ) \times_{P_{h,d-h}(\OC_v)} GL_d(\OC_v)$ is also the intermediate
%extension of  $HT_{\overline{1_h}}(\chi_v,\st_h(\chi_v))  
%\times_{P_{h,d-h}(\OC_v)} GL_d(\OC_v)$.
%\end{proof}

By adjunction property, the map
\addtocounter{thm}{1}
\begin{multline} \label{eq-map1}
j_!^{=h+\delta} HT(\chi_v,\st_h(\chi_v \{ \frac{-\delta}{2} \} ) \times \speh_{\delta}
(\chi_v \{ h/2 \} )) \otimes \Xi^{\delta/2} \\
\longrightarrow j_!^{=h+\delta-1} HT(\chi_v,\st_h(\chi_v \{ \frac{1-\delta}{2} \} ) 
\times  \speh_{\delta-1}(\chi_v\{ h/2 \} )) \otimes \Xi^{\frac{\delta-1}{2}}
\end{multline}
is given by 
\addtocounter{thm}{1}
\begin{multline}\label{eq-map2}
HT(\chi_v,\st_h(\chi_v \{ \frac{-\delta}{2} \} ) \times \speh_{\delta}(\chi_v\{ h/2 \} )) 
\otimes \Xi^{\delta/2} \longrightarrow \\ j^{=h+\delta,*} (
\lexp p i^{h+\delta,!}  ( 
j_!^{=h+\delta-1} HT(\chi_v,\st_h(\chi_v \{ \frac{1-\delta}{2} \} ) 
\times \speh_{\delta-1}(\chi_v\{ h/2 \} )) \otimes \Xi^{\frac{\delta-1}{2}}))
\end{multline}
To compute this last term we use the resolution (\ref{eq-resolution0}).
Precisely denote by 
$\HC:=HT(\chi_v,\st_h(\chi \{ \frac{1-\delta}{2} \} ) \times 
\speh_{\delta-1}(\chi_v\{ h/2 \} )) \otimes \Xi^{\frac{\delta-1}{2}}$,
and write the previous resolution as follows
$$0 \rightarrow K \longrightarrow  j_!^{=h+\delta} \HC' \longrightarrow Q \rightarrow 0,$$
$$0 \rightarrow Q \longrightarrow  j_!^{=h+\delta-1} \HC 
\longrightarrow  \lexp p j_{!*}^{=h+\delta-1} \HC \rightarrow 0,$$
with 
$$\HC':=HT \Bigl ( \chi_v,  \st_h(\chi_v \{ \frac{1-\delta}{2} \} ) 
\times \bigl ( \speh_{\delta-1}(\chi_v \{ -1/2 \})
 \times \chi_v \{\frac{\delta-1}{2} \}  \bigr ) \{ h/2 \}  \Bigr ) \otimes \Xi^{\delta/2}.$$ 
As the support of $K$ is contained in $\sh^{\geq h+\delta+1}_{I,\bar s_v}$ then
$\lexp p i^{h+\delta,!} K=K$ and
$j^{=h+\delta,*} (\lexp p i^{h+\delta,!} K)$ is zero. Moreover
$\lexp p i^{h+\delta,!} ( \lexp p j_{!*}^{=h+\delta-1} \HC)$ is zero by
construction of the intermediate extension. We then deduce that
\addtocounter{thm}{1}
\begin{multline} \label{eq-map3}
 j^{=h+\delta,*} (\lexp p i^{h+\delta,!} ( j_!^{=h+\delta-1} HT(\chi_v,\st_t(\chi_v 
 \{ \frac{1-\delta}{2} \} ) \times 
\speh_{\delta-1}(\chi_v\{ h/2 \} )) \otimes \Xi^{\frac{\delta-1}{2}})) \\ \simeq 
HT \Bigl ( \chi_v,  \st_h(\chi_v \{ \frac{1-\delta}{2} \} ) \\ \times \bigl 
( \speh_{\delta-1}(\chi_v \{ -1/2 \})
 \times \chi_v \{\frac{\delta-1}{2} \}  \bigr ) \{ h/2 \}  \Bigr ) \otimes \Xi^{\delta/2} 
\end{multline}
%
%When $\pi_{v,-1}$ is a character $\chi_v$, the previous isomorphism is easily 
%seen. Indeed 
%
%
%from (\ref{eq-resolution0}) and before inducing from
%$\sh^{=h}_{I,\bar s_v,\overline{1_h}}$ to $\sh^{=h}_{I,\bar s_v}$, we just need to understand  
%$\lexp p i^{h+1,!} \lexp p j^{=h}_{\overline{1_h},!*} HT(\chi_v,\Pi_h)$
%
%
% knowning,
%cf. the proof of the lemma \ref{lem-ext0}, that, 
%as $\sh^{\geq h}_{I,\bar s_v,\overline{1_h}}$ is smooth,
%$\lexp p j^{=h}_{\overline{1_h},*} HT(\chi_v,\Pi_h)$ is isomorphic, up to
%twist the action of the fundamental group by a character, to
%$ (\overline \Zm_l)_{|\sh^{\geq h}_{I,\bar s_v,\overline{1_h}}} \otimes \Pi_h.$
%
\marque \emph{Fact}.
In particular, up to homothety, the map (\ref{eq-map3}), and so those of (\ref{eq-map2}), is unique.
Finally as the maps of (\ref{eq-resolution0}) are strict, the given maps (\ref{eq-map1}) are uniquely
determined, that is, if we forget the infinitesimal parts, these maps are independent of the 
chosen $t$ in (\ref{eq-resolution0}).

For every $1 \leq h \leq d$, let denote by $i(h)$ the smallest index $i$ such that 
$H^i(\sh_{I,\bar s_v},\lexp p j^{=h}_{!*} HT(\chi_v ,
\st_h(\chi_v)))_{\mathfrak m}$
has non trivial torsion: if it does not exist then we set $i(h)=+\oo$. By duality, as 
$\lexp p j_{!*}=\lexp {p+} j_{!*}$ for Harris-Taylor local systems associated
to characters, note that when $i(h)$ is finite then $i(h) \leq 0$. 
Suppose by absurdity there exists
$h$ with $i(h)$ finite and denote $h_0$ the biggest such $h$.

\begin{lemma} \label{lem-ih}
For $1 \leq h \leq h_0$ then $i(h)=h-h_0$.
\end{lemma}

\rem A similar result is proved in \cite{boyer-imj} when the
level is maximal at $v$.

\begin{proof}
a) We first prove that for every $h_0 < h \leq d$, the cohomology 
groups of  $j^{=h}_! HT(\chi_v,\Pi_h)$ are torsion free. Consider the
following strict filtration in the category of free perverse sheaves
\begin{multline} \label{eq-fil-j}
(0)=\Fil^{-1-d}(\chi_v ,h) \harrow \Fil^{-d}(\chi_v ,h) \harrow \cdots \\
\harrow \Fil^{-h}(\chi_v ,h)=j^{=h}_{!} HT(\chi_v ,\Pi_h)
\end{multline}
where the symbol $\harrow$ means a strict monomorphism, 
with graded parts 
$$\gr^{-k}(\chi_v,h) \simeq \lexp p j^{=k}_{!*} 
HT(\chi_v,\Pi_h \{\frac{h-k}{2} \} \otimes \st_{k-h}(\chi_v\{h/2 \} ))(\frac{h-k}{2}).$$ 
Over $\overline \Qm_l$, the result is proved in
\cite{boyer-invent2} \S 4.3. From \cite{boyer-torsion} such a filtration
can be constructed over $\overline \Zm_l$ up to the fact that the graduate parts
are only known to verify
\begin{multline*}
\lexp p j^{=k}_{!*} 
HT(\chi_v,\Pi_h \{\frac{h-k}{2} \} \otimes \st_{k-h}(\chi_v \{h/2 \} ))(\frac{h-k}{2})
\htarrow \gr^{-k}(\chi_v,h) \\
 \htarrow  \lexp {p+} j^{=k}_{!*} 
HT(\chi_v,\Pi_h \{\frac{h-k}{2} \} \otimes \st_{k-h}(\chi_v\{h/2 \} ))(\frac{h-k}{2}),
\end{multline*}
and we can conclude thanks to (\ref{eq-ext0}). The associated
spectral sequence localized at $\mathfrak m$, is then concentrated in middle degree and torsion free which gives the claim.

b) Before watching the cases $h \leq h_0$, note that
the spectral sequence associated to (\ref{eq-resolution0}) for $h=h_0+1$, 
has all its $E_1$ terms torsion free
and degenerates at its $E_2$ terms. As by hypothesis the aims of this spectral sequence is free
and equals to only one $E_2$ terms, we deduce that all the maps
\addtocounter{thm}{1}
\begin{multline} \label{eq-map1-coho}
H^0 \bigl (\sh_{I,\bar s_v},j_!^{=h+\delta} HT_\xi(\chi_v,\st_h(\chi_v \{ \frac{-\delta}{2} \} ) \times \speh_{\delta}
(\chi_v\{ h/2 \} )) \otimes \Xi^{\delta/2} \bigr )_{\mathfrak m} \\
\longrightarrow \\ H^0 \bigl (\sh_{I,\bar s_v},
j_!^{=h+\delta-1} HT_\xi(\chi_v,\st_h(\chi_v \{ \frac{1-\delta}{2} \} ) \\ \times 
\speh_{\delta-1}(\chi_v\{ h/2 \} )) \otimes \Xi^{\frac{\delta-1}{2}} \bigr )_{\mathfrak m}
\end{multline}
are saturated, i.e. their cokernel are free $\overline \Zm_l$-modules. Then from the previous fact stressed after (\ref{eq-map3}), this property
remains true when we consider the associated spectral sequence for 
$1 \leq h' \leq h_0$.

c) Consider now $h=h_0$ and the spectral sequence associated to 
(\ref{eq-resolution0}) where
\addtocounter{thm}{1}
\begin{multline} \label{eq-E2pq}
E_2^{p,q}=H^{p+2q}(\sh_{I,\bar s_v}, j_!^{=h+q} \\ 
HT_\xi(\chi_v,\st_h(\chi_v (-q/2)) \times \speh_q(\chi_v \{ h/2 \} )) 
\otimes \Xi^{\frac{q}{2}})_{\mathfrak m}
\end{multline}
By definition of $h_0$, we know that  some of
the $E_\oo^{p,-p}$ should have a non trivial torsion subspace.
We saw that 
\begin{itemize}
\item the contributions from the deeper strata are torsion free and

\item $H^i(\sh_{I,\bar s_v},j^{=h_0}_! HT_\xi(\chi_v,\Pi_{h_0}))_{\mathfrak m}$
are zero for $i<0$ and is torsion free for $i=0$, whatever is $\Pi_{h_0}$.

\item Then there should exist
a non strict map $d_1^{p,q}$. But, we have just seen that it
can not be maps between deeper strata.

\item Finally, using the previous points, the only possibility is that the 
cokernel of
\addtocounter{thm}{1}
\begin{multline} \label{eq-map1-coho2}
H^0 \bigl (\sh_{I,\bar s_v},j_!^{=h_0+1} HT_\xi(\chi_v,\st_{h_0}(\chi_v 
\{ \frac{-1}{2} \} ) 
\times \chi_v \{ h_0/2 \} )) \otimes \Xi^{1/2} \bigr )_{\mathfrak m} \\
\longrightarrow \\ H^0 \bigl (\sh_{I,\bar s_v},
j_!^{=h_0} HT_\xi(\chi_v,\st_{h_0}(\chi_v)) \bigr )_{\mathfrak m}
\end{multline}
has a non trivial torsion subspace. 
\end{itemize}
In particular we have $i(h_0)=0$.

d) Finally using the fact 2.18 and the previous points, for any $1 \leq h \leq h_0$, 
in the spectral sequence (\ref{eq-E2pq})
\begin{itemize}
\item by point a), $E_2^{p,q}$ is torsion free for $q \geq h_0-h+1$ and so it
is zero if $p+2q \neq 0$;

\item by affiness of the open strata, cf. \cite{boyer-imj} theorem 1.8,
$E_2^{p,q}$ is zero for $p+2q<0$ and torsion free for $p+2q=0$;

\item by point b), the maps $d_2^{p,q}$ are saturated for $q \geq h_0-h+2$;

\item by point c), $d_2^{-2(h_0-h+1),h_0-h+1}$ has a cokernel with a non trivial 
torsion subspace.

\item Moreover, over $\overline \Qm_l$, the spectral sequence degenerates
at $E_3$ and $E_3^{p,q}=0$ if $(p,q) \neq (0,0)$.
\end{itemize}
We then deduce that 
$H^i(\sh_{I,\bar s_v},\lexp p j^{=h}_{!*} HT_\xi(\chi_v,\Pi_h))_{\mathfrak m}$
is zero for $i < h-h_0$ and for $i=h-h_0$ it has a non trivial torsion subspace.
\end{proof}

%
%\rem Proposition \ref{prop-dec-pervers} implies that for any $\pi_v$ of type $\varrho$ and 
%for every $t$,
%the torsion of $H^i(\sh_{I,\bar s_v},\lexp p j^{=tg}_{!*} HT(\pi_v,\Pi_t))_{\mathfrak m}$
%is trivial for any $i \leq 1-t_0$ which gives us the corresponding informations about the 
%associated map (\ref{eq-map1-coho}).
%

Consider now the filtration of stratification of $\Psi_\varrho$ constructed using the
adjunction morphisms $j^{=h}_! j^{=h,*}$ as in \cite{boyer-torsion}
\addtocounter{thm}{1}
\begin{equation} \label{eq-fil-psi}
\Fil^1_!(\Psi_\varrho) \harrow \Fil^{2}_!(\Psi_\varrho)
\harrow \cdots \harrow \Fil^{d}_!(\Psi_\varrho)
\end{equation}
where $\Fil^{h}_!(\Psi_\varrho)$ is the saturated image of $j^{=h}_!j^{=h,*} \Psi_\varrho
\longrightarrow \Psi_\varrho$.

\rem Recall that the filtration
$\Fill^\bullet$ is a refinement of $\Fil_!^\bullet$ as one can see it in the next proposition.

For our fixed $\chi_v$, let denote 
$\Fil^1_{!,\chi_v}(\Psi) \harrow \Fil^1_!(\Psi_\varrho)$ such that
$\Fil^1_{!,\chi_v}(\Psi) \otimes_{\overline \Zm_l} \overline \Qm_l \simeq \Fil^1_!(\Psi_{\chi_v})$
where $\Psi_{\chi_v}$ is the direct factor of $\Psi \otimes_{\overline \Zm_l} \overline \Qm_l$
associated to $\chi_v$, cf. \cite{boyer-torsion}. 

\begin{prop} (cf. \cite{boyer-duke} 3.3.5)
We have the following resolution of $\gr^h_{!,\chi_v}(\Psi)$
\addtocounter{thm}{1}
\begin{multline} \label{eq-resolution-psi2}
0 \rightarrow j^{=d}_! HT(\chi_v,LT_{h,d}(\chi_v)) \otimes L_g(\chi_v(\frac{d-h}{2}))
\longrightarrow \\
 j^{=d-1}_!HT(\chi_v,LT_{h,d-1}(\chi_v)) \otimes L_g(\chi_v(\frac{d-h-1}{2}) )
\longrightarrow \\ \cdots \longrightarrow j^{=h}_! HT(\chi_v,\st_h(\chi_v)) \otimes 
\Lm(\chi_v) \longrightarrow \gr^{h}_{!,\chi_v}(\Psi) \rightarrow 0,
\end{multline}
where
\begin{itemize}
\item $LT_{h,h+\delta}(\chi_v) \hookrightarrow \st_h(\chi_v \{ -\delta/2 \}) \times \speh_{\delta}(\chi_v \{ h/2 \}),$
is the only irreducible sub-space of this induced representation,

\item and $\Lm$ is the local Langlands correspondence composed
by contragredient.
\end{itemize}
%
%
%\addtocounter{thm}{1}
%\begin{multline} \label{eq-resolution-psi}
%0 \rightarrow j^{=d}_! HT(\chi_v,\speh_d(\chi_v)) \otimes \Lm(\chi_v(\frac{d-1}{2}))
%\longrightarrow \\
% j^{=d-1}_!HT(\chi_v,\speh_{d-1}(\chi_v)) \otimes \Lm(\chi_v(\frac{d-2}{2}) )
%\longrightarrow \\ \cdots \longrightarrow j^{=1}_! HT(\chi_v,\chi_v) \otimes 
%\Lm(\chi_v) \longrightarrow \Fil^1_{!,\chi_v}(\Psi_v) \rightarrow 0,
%\end{multline}
%where we recall that $\Lm$ is the local Langlands correspondence composed
%by contragredient.
\end{prop}

\rems 
\begin{itemize}
\item As explained after proposition \ref{prop-resolution1},
it amounts to describe the germs of the $\overline \Zm_l$-sheaf cohomology of 
$\gr^h_{!,\chi_v}(\Psi_{v,\xi})$.
Over $\overline \Qm_l$, the resolution (\ref{eq-resolution-psi2}) is then proved in
\cite{boyer-invent2}.

\item Over $\overline \Zm_l$,
it is proved in full generality in \cite{boyer-duke} for every irreducible
supercuspidal representation $\pi_v$ in place of $\chi_v$. It amounts to prove
that the germs of the sheaf cohomology of 
$\gr^1_{!,\chi_v}(\Psi_{v,\xi})$ are free. The case of a character is however much 
more simple. Indeed consider then the torsion part of the cokernel of one 
of these maps.
Note that, thanks to (\ref{eq-ext0}), such a cokernel 
must have non trivial invariants under the action 
the Iwahori sub-group at $v$. We then work at Iwahori level at $v$.
As said above, it amounts to understand the germs of the 
$\overline \Zm_l$-sheaf cohomology of $\gr^h_{!,\chi_v}(\Psi)$ 
which are described, cf. \cite{fargues-annexe},
by the cohomology of the Lubin-Tate tower. By the comparison theorem of 
Faltings-Fargues, cf. \cite{fargues-faltings}, one is reduced to compute the cohomology of the Drinfeld tower in Iwahori level which is already done 
in \cite{s-s}: we then note that there are all free $\overline \Zm_l$-modules.
\end{itemize}
We can then apply the previous arguments a)-d) above, 
for $h \leq h_0$ (resp. $h > h_0$) 
the torsion of $H^i(\sh_{I,\bar s_v},\gr^{h}_{!,\chi_v}(\Psi_{v,\xi}))_{\mathfrak m}$
is trivial for any $i \leq h-h_0$ (resp. for all $i$) 
and the free parts are concentrated for $i=0$. 
Using then the spectral sequence associated to the previous filtration, 
we can then conclude that $H^{1-t_0}(\sh_{I,\bar s_v},\Psi_{v,\xi})_{\mathfrak m}$ 
would have non trivial torsion 
which is false as $\mathfrak m$ is supposed to be KHT-free.
\end{proof}

\section{Local behavior of monodromy over $\overline \Fm_l$}
\label{para-local}

Recall that $\varrho$ is a fixed $\overline \Fm_l$-character and $\Psi_\varrho$ is the
associated direct factor of $\Psi_v$. Over $\overline \Qm_l$, the monodromy operator
define a nilpotent morphism 
$N_{\varrho,\overline \Qm_l}: \Psi_\varrho \otimes_{\overline \Zm_l} \overline \Qm_l
\longrightarrow \Psi_\varrho  \otimes_{\overline \Zm_l} \overline \Qm_l$
compatible with the filtration $\Fil^\bullet_!(\Psi_\varrho)$ in the sense that
$\Fil^h_!(\Psi_\varrho)  \otimes_{\overline \Zm_l} \overline \Qm_l$ coincides
with the kernel of $N_{\varrho,\overline \Qm_l}^h$. The aim of this section is to
construct a $\overline \Zm_l$-version $N_\varrho$ of $N_{\varrho,\overline \Qm_l}$ such that
$\Fil^h_!(\Psi_\varrho)  \otimes_{\overline \Zm_l} \overline \Fm_l$ coincides with
the kernel of $N_\varrho^h  \otimes_{\overline \Zm_l} \overline \Fm_l$.

\bigskip

\noindent \textit{First step}: consider 
$$0 \rightarrow \Fil^{1}_!(\Psi_\varrho) \longrightarrow \Psi_\varrho \longrightarrow
\coFil^1_!(\Psi_\varrho) \rightarrow 0,$$
and the following long exact sequence
\begin{multline*}
0 \rightarrow \hom(\coFil^1_!(\Psi_\varrho),\Psi_\varrho) \longrightarrow
\hom(\Psi_\varrho,\Psi_\varrho) \\ \longrightarrow \hom(\Fil^1_!(\Psi_\varrho),\Psi_\varrho)
\longrightarrow \cdots
\end{multline*}
where $\hom$ is taken
in the category of equivariant Hecke perverse sheaves.\footnote{We do not ask the map to be Galois equivariant as $N_{\varrho,\overline \Qm_l}$ is not.} Note that
$N_{\varrho,\overline \Qm_l} \in \hom(\Psi_\varrho,\Psi_\varrho) \otimes_{\overline \Zm_l}
\overline \Qm_l$ comes from 
$\hom(\coFil^1_!(\Psi_\varrho),\Psi_\varrho)\otimes_{\overline \Zm_l}
\overline \Qm_l$, so that we focus on  $\hom(\coFil^1_!(\Psi_\varrho),\Psi_\varrho)$.
From
$$0 \rightarrow \gr^2_!(\Psi_\varrho) \longrightarrow \coFil^1_!(\Psi_\varrho) 
\longrightarrow \coFil^2_!(\Psi_\varrho) \rightarrow 0,$$
we obtain
\begin{multline*}
0 \rightarrow \hom(\coFil_!^2(\Psi_\varrho),\Psi_\varrho) \longrightarrow
\hom(\coFil_!^1,\Psi_\varrho) \longrightarrow \\
\hom(\gr_!^2(\Psi_\varrho),\Psi_\varrho) \longrightarrow \ext^1(\coFil_!^2(\Psi_\varrho),\Psi_\varrho)) \longrightarrow \cdots
\end{multline*}
Note then that $N_{\varrho,\overline \Qm_l} \in \hom(\coFil_!^1,\Psi_\varrho)
\otimes_{\overline \Zm_l} \overline \Qm_l$ does not belong to the image of
$\hom(\coFil_!^2(\Psi_\varrho),\Psi_\varrho)\otimes_{\overline \Zm_l}
\overline \Qm_l$.

\begin{lemma}
The $\overline \Zm_l$-module $\ext^1(\coFil_!^2(\Psi_\varrho),\Psi_\varrho))$
is torsion free.
\end{lemma}

\begin{proof}
Let $\alpha \in \ext^1(\coFil_!^2(\Psi_\varrho),\Psi_\varrho))$ which is killed
by some power of $l$ and let
$$0 \rightarrow \Psi_\varrho \longrightarrow P \longrightarrow
\coFil_!^2(\Psi_\varrho) \rightarrow 0,$$
be the extension defined by $\alpha$. Applying 
$\otimes^\Lm_{\overline \Zm_l} \overline \Qm_l$, this short exact sequence split so that
$P$ can be written
$$0 \rightarrow \widetilde \coFil_!^2(\Psi_\varrho) \longrightarrow P \longrightarrow
\widetilde \Psi_\varrho  \rightarrow 0,$$
so that composing through $P$ we obtain
$$0 \rightarrow \widetilde \coFil_!^2(\Psi_\varrho) \longrightarrow \coFil_!^2(\Psi_\varrho)
\longrightarrow T \rightarrow 0,$$
$$0 \rightarrow \Psi_\varrho \longrightarrow 
\widetilde \Psi_\varrho \longrightarrow T \rightarrow 0,$$
for the same torsion perverse sheaf $T$ appearing as cokernel of the two previous maps.
By tensoring with $\otimes_{\overline \Zm_l} \overline \Fm_l$, we then obtain
\begin{multline*}
0 \rightarrow \lexp ph^{-1} (T \otimes_{\overline \Zm_l} \overline \Fm_l) \longrightarrow
\widetilde \coFil_!^2(\Psi_\varrho) \otimes_{\overline \Zm_l} \overline \Fm_l \\
 \longrightarrow \coFil_!^2(\Psi_\varrho) \otimes_{\overline \Zm_l} \overline \Fm_l
\longrightarrow \lexp p h^0(T\otimes_{\overline \Zm_l} \overline \Fm_l) \rightarrow 0,
\end{multline*}
and
\begin{multline*}
0 \rightarrow \lexp ph^{-1} (T \otimes_{\overline \Zm_l} \overline \Fm_l) \longrightarrow
\Psi_\varrho \otimes_{\overline \Zm_l} \overline \Fm_l \\
 \longrightarrow \widetilde \Psi_\varrho \otimes_{\overline \Zm_l} \overline \Fm_l
\longrightarrow \lexp p h^0(T\otimes_{\overline \Zm_l} \overline \Fm_l) \rightarrow 0,
\end{multline*}
where we denote by $\lexp p h^\bullet K$ the $p$-perverse cohomology complexes
of $K$. As by lemma \ref{lem-socle}, 
the socle of $\Psi_\varrho \otimes_{\overline \Zm_l} \overline \Fm_l$
is, up to multiplicity, 
$j^{=d}_{!*} HT_{\overline \Fm_l}(\varrho,\st_d(\varrho))(\frac{d-1}{2})$ then
it has to be a constituant of $\lexp ph^{-1} (T \otimes_{\overline \Zm_l} \overline \Fm_l)$.
But, using that the order of $q_v$ modulo $l$ is strictly greater than $d$,
$j^{=d}_{!*} HT_{\overline \Fm_l}(\varrho,\st_d(\varrho))(\frac{d-1}{2})$ is not a constituant
of $ \coFil_!^2(\Psi_\varrho) \otimes_{\overline \Zm_l} \overline \Fm_l$.
We then deduce that $T$ is forced to be zero, which means that the extension $P$
is split, i.e. $\alpha \in \ext^1(\coFil^2_!(\Psi_\varrho),\Psi_\varrho)$
is zero.
\end{proof}

We are then led to study
$$\hom(\gr_!^2(\Psi_\varrho),\Psi_\varrho) \simeq
\hom(\gr_!^2(\Psi_\varrho), \Fil_*^1(\gr_!^1(\Psi_\varrho)))$$
where
$$0 \rightarrow \Fil_*^1(\gr_!^1(\Psi_\varrho)) \longrightarrow \Fil_!^1(\Psi_\varrho) 
\longrightarrow \coFil_*^1(\Fil_!^1(\Psi_\varrho) \rightarrow 0.$$
Note that, up to an unramified Galois twist, 
$\gr_!^2(\Psi_\varrho)\otimes_{\overline \Zm_l} \overline \Qm_l \simeq 
\Fil_*^1(\gr_1^1(\Psi_\varrho)) \otimes_{\overline \Zm_l} \overline \Qm_l$ and the
cosocle of $\gr_!^2(\Psi_\varrho)\otimes_{\overline \Zm_l} \overline \Fm_l$ (resp. $\Fil_*^1(\gr_1^1(\Psi_\varrho))\otimes_{\overline \Zm_l} \overline \Fm_l$)
is some multiple, depending on the lattice of $j^{=2,*} \gr_!^2(\Psi_\varrho)$
(resp. $j^{=2,*} \Fil_*^1(\gr_1^1(\Psi_\varrho))$, of
$j^{=2}_{!*} HT(\varrho,\st_2(\varrho))$. 
Consider then
$$0 \rightarrow P_- \longrightarrow \Fil_*^1(\gr_1^1(\Psi_\varrho)) \longrightarrow
\lexp p j^{=2}_{!*} j^{=2,*} \Fil_*^1(\gr_1^1(\Psi_\varrho)) \rightarrow 0,$$
where $P_-$ has support in $\sh^{\geq 3}_{v}$. Then using as before long exact 
sequence, we note that
$$\hom(\gr_!^2(\Psi_\varrho), \Fil_*^1(\gr_1^1(\Psi_\varrho))) \simeq
\hom(\lexp p j^{=2}_{!*} j^{=2,*}\gr_!^2(\Psi_\varrho),\lexp p j^{=2}_{!*} j^{=2,*}
\Fil_*^1(\gr_1^1(\Psi_\varrho))).$$ 
In particular if the two local systems
$\lexp p j^{=2}_{!*} j^{=2,*}\gr_!^2(\Psi_\varrho)$ and $\lexp p j^{=2}_{!*} j^{=2,*}
\Fil_*^1(\gr_1^1(\Psi_\varrho))$ were isomorphic, then there exist a element
in $\hom(\gr_!^2(\Psi_\varrho), \Fil_*^1(\gr_1^1(\Psi_\varrho)))$ which is an isomorphism.
This element then gives us a $\overline \Zm_l$-morphism 
$N_v \in \hom( \Psi_\varrho,\Psi_\varrho)$ so that $\Fil^1_*(\gr^1_!(\Psi_\varrho))$
is in the image.

\bigskip

\noindent \textit{Second step}: we want to prove that the local systems
$j^{=2,*}\gr_!^2(\Psi_\varrho)$ and $j^{=2,*}
\Fil_*^1(\gr_!^1(\Psi_\varrho))$ are isomorphic. Consider first the following situation:
let $\LC_k$ and $\LC_{k+1}$ be $\overline \Zm_l$-local systems on a scheme $X$ 
such that:
\begin{itemize}
\item $\gr_{k+1,\overline \Qm_l}:=(\LC_{k+1}/\LC_k) \otimes_{\overline \Zm_l} \overline \Qm_l$ is irreducible and we introduce
$$\xymatrix{
\gr_{k+1} \ar@{^{(}-->}[r] \ar@{^{(}-->}[d] & \gr_{k+1,\overline \Qm_l} \ar@{^{(}->}[d] \\
\LC_{k+1} \ar@{^{(}->}[r] & \LC_{k+1} \otimes_{\overline \Zm_l} \overline \Qm_l.
}$$
We moreover suppose that the $\gr_{k+1} \otimes_{\overline \Zm_l} \overline \Fm_l$
are also irreducible so the various stable lattices of $\gr_{k+1}$ are homothetic.

\item $\LC_{k+1} \otimes_{\overline \Zm_l} \overline \Qm_l \simeq
\LC_k \otimes_{\overline \Zm_l} \overline \Qm_l \oplus \gr_{k+1} \otimes_{\overline \Zm_l} \overline \Qm_l$.
\end{itemize}
We then have
$$0 \rightarrow \LC_k \oplus \gr_{k+1} \longrightarrow \LC_{k+1} 
\longrightarrow T \rightarrow 0,$$
where $T$ is torsion and can be viewed as a quotient
$$\LC_k \hookrightarrow \LC'_k \twoheadrightarrow T, \quad 
\gr_{k+1} \hookrightarrow \gr'_{k+1} \twoheadrightarrow T,$$
with, cf. also the exchange basic step of \S \ref{para-notations}
$$\LC_k \hookrightarrow \LC_{k+1} \twoheadrightarrow \gr'_{k+1}, \qquad
\gr_{k+1} \hookrightarrow \LC_{k+1} \twoheadrightarrow \LC'_k.$$
As $\gr_{k+1} \otimes_{\overline \Zm_l} \overline \Qm_l$ is irreducible, then
$\gr_{k+1} \hookrightarrow \gr'_{k+1}$ is given by multiplication by $l^\delta$ and the
extension is characterized by this $\delta$.

Consider then the $\overline \Zm_l$-local system $\LC:=j^{=1,*} \Psi_\varrho$
and recall that 
$$\LC \otimes_{\overline \Zm_l} \overline \Qm_l \simeq 
\bigoplus_{i=1}^r HT_{\overline \Qm_l}(\chi_{v,i},\chi_{v,i}),$$
where we fix any numbering of $\cusp(\varrho)=\{\chi_{v,1},\cdots,\chi_{v,r}\}$. 
For $k=1,\cdots,r$, we introduce
$$\xymatrix{
\LC_k \ar@{^{(}-->}[r] \ar@{^{(}-->}[d] & \bigoplus_{i=1}^k HT(\chi_{v,i},\chi_{v,i}) 
\ar@{^{(}->}[d] \\
\LC \ar@{^{(}->}[r] & \LC \otimes_{\overline \Zm_l} \overline \Qm_l.
}$$
Let denote by $T_{k+1}$ the torsion local system such that
$$0 \rightarrow \LC_k \oplus \gr_{k+1} \longrightarrow \LC_{k+1} \longrightarrow T_{k+1}
\rightarrow 0,$$
where $\gr_{k+1}:=\LC_{k+1}/\LC_k$, as above. We can apply the previous remark
and denote by $\delta_k$ the power of $l$ which define the homothety
$\gr_{k+1} \hookrightarrow \gr'_{k+1} \twoheadrightarrow T_{k+1}$. The set of
$\delta_k$ for $k=1,\cdots,r$ is then a numerical data to characterize $\LC$ inside
$j^{=1,*} \Psi_\varrho \otimes_{\overline \Zm_l} \overline \Qm_l$. 

(i) From the main result of \cite{boyer-duke}, 
cf. its introduction paragraph, $\Fil_*^1(\gr^1_!(\Psi_\varrho))$ is obtained as follows
$$0 \rightarrow \Fil_*^1(\gr^1_!(\Psi_\varrho)) \longrightarrow 
j^{=1}_{\neq 1,!} j^{=1,*}_{\neq 1} \Psi_\varrho \longrightarrow  
\lexp p j^{=1}_{\neq 1,!*} j^{=1,*}_{\neq 1} \Psi_\varrho \rightarrow 0,$$
so that $j^{=2,*} \Fil_*^1(\gr^1_!(\Psi_\varrho)) \simeq \lexp p h^{-1} i^{2,*}
\lexp p j^{=1}_{\neq 1,!*} j^{=1,*}_{\neq 1} \Psi_\varrho$.
With the previous notations, we have
$$0 \rightarrow \lexp p j^{=1}_{\neq 1,!*} \LC_k \longrightarrow 
\lexp p j^{=1}_{\neq 1,!*} \LC_{k+1} \longrightarrow \lexp p j^{=1}_{\neq 1,!*} T_{k+1} 
\rightarrow 0,$$
from which we obtain the following description of $j^{=2,*} \Fil_*^1(\gr^1_!(\Psi_\varrho))$:
\begin{itemize}
\item there exists local systems $\LC^+_k$ for $k=1,\cdots,r$
so that $\LC^+_k \otimes_{\overline \Zm_l} \overline \Qm_l \simeq \bigoplus_{i=1}^k
HT_{\overline \Qm_l}(\chi_{v,i},\st_2(\chi_{v,i})(-1/2)$;

\item with $\gr^+_{k+1}$ defined, as before, with
$$0 \rightarrow \LC^+_k \oplus \gr^+_{k+1} \longrightarrow \LC^+_{k+1}
\longrightarrow T_{k+1} \rightarrow 0,$$
where $T_{k+1}$ is killed by $l^{\delta_{k+1}}$.
\end{itemize}
Finally $j^{=2,*} \Fil_*^1(\gr^1_!(\Psi_\varrho))$ is described with the same numerical
data $\{ \delta_k: k=1,\cdots,r \}$ as $j^{=1,*} \Psi_\varrho$.

(ii) The same arguments apply with $j^{=1}_{\neq 1,*}j^{=1,*}_{\neq 1} \Psi_\varrho$
so that the local system $\widetilde \LC:=\lexp {p+} h^1i^{2,*} \lexp {p+} j^{=1}_{\neq 1,!*}
j^{=1,*}_{\neq 1} \Psi_\varrho$ is also characterized by the same numerical data
$\{ \delta_k: k=1,\cdots,r \}$ except that $\widetilde \LC$ can not directly be
identified to $j^{=2,*} \gr^2_!(\Psi_\varrho)$. Indeed we are interested in the lattice of
$\bigoplus_{\chi_v \in \cusp(\varrho)} j^{=2}_{!*} 
HT_{\overline \Qm_l}(\chi_v,\st_2(\chi_v))(-1/2)$ given by 
$\Psi_\varrho/\Fil^1_!(\Psi_\varrho)$.  But by now the
previous lattice of $P_1:=\lexp p j^{=2}_{!*} \widetilde \LC$ described by 
$\{ \delta_k: k=1,\cdots r\}$ is obtained using a filtration where $P_1$ appears as the
socle of the perverse sheaf $Q$ defined as follows:
$$0 \rightarrow \lexp p j^{=1}_{!*} j^{=1,*} \coFil^1_*(\Psi_\varrho) \longrightarrow
\coFil^1_*(\Psi_\varrho) \longrightarrow Q \rightarrow 0.$$
As explained in \S \ref{para-notations}, we have to use basic exchange
steps as many times as needed to move $P_1$ until it appears 
as the cosocle of $\Fil^2_!(\Psi_\varrho) \hookrightarrow \Psi_\varrho$,
cf. the discussion before lemma \ref{lem-socle}.

Note then that all the perverse sheaves which are exchanged with $P_1$ during
this process, are lattice of
$j^{=h}_{!*} HT_{\overline \Qm_l}(\chi_v,\st_h(\chi_v))(\frac{1-h+\delta}{2})$ with
$h \geq 3$. As explained in the remark after the definition of the exchange basic step,
as $\lexp p j^{=2}_{!*} HT(\chi_v,\st_2(\chi_v)) \simeq
\lexp {p+} j^{=2}_{!*} HT(\chi_v,\st_2(\chi_v))$, for all these exchange, we have
$T=0$ and $P_1$ remains unchanged during all the basic exchange steps.

\medskip

\noindent \textit{Third step}: at this stage we constructed a $\overline \Zm_l$-version
of the monodromy operator 
$$\xymatrix{
N_\varrho: \Psi_\varrho \ar[rr] \ar[dr] & & \Psi_\varrho \\
& \Psi_\varrho/\Fil^1_!(\Psi_\varrho) \ar[ur]
}$$
such that the kernel of 
$$N_\varrho \otimes_{\overline \Zm_l} \overline \Fm_l: \bigl ( \Psi_\varrho/\Fil^1_!(\Psi_\varrho)
\bigr ) \otimes_{\overline \Zm_l} \overline \Fm_l \longrightarrow \Psi_\varrho 
\otimes_{\overline \Zm_l} \overline \Fm_l,$$
does not contain any irreducible subquotient of 
$\gr^2_!(\Psi_\varrho) \otimes_{\overline \Zm_l} \overline \Fm_l$.

Recall that we suppose the order of $q_v$ modulo $l$ to be $>d$, so that 
the irreducible constituants of 
$\gr^2_!(\Psi_\varrho) \otimes_{\overline \Zm_l} \overline \Fm_l$
are disjoint from that of $\bigl ( \Psi_\varrho/\Fil2_!(\Psi_\varrho) \bigr) 
\otimes_{\overline \Zm_l} \overline \Fm_l$.

Moreover, arguing as in lemma \ref{lem-socle}, we see that the socle of
$\bigl (\Psi_\varrho/\Fil^1_!(\Psi_\varrho) \bigr ) \otimes_{\overline \Zm_l} \overline \Fm_l$ is up to
multiplicity $j^{=d}_{!*} HT_{\overline \Fm_l}(\varrho,\st_d(\varrho))(\frac{d-3}{2})$
which is a constituant
of $\gr^2_!(\Psi_\varrho) \otimes_{\overline \Zm_l} \overline \Fm_l$.

From the previous facts, we then deduce that the kernel of
$N_\varrho \otimes_{\overline \Zm_l} \overline \Fm_l: \Psi_\varrho \otimes_{\overline \Zm_l} \overline \Fm_l \longrightarrow \Psi_\varrho \otimes_{\overline \Zm_l} \overline \Fm_l$
is reduced to $\Fil^1_!(\Psi_\varrho) \otimes_{\overline \Zm_l} \overline \Fm_l$ and is
so the modulo $l$ reduction of the kernel of $N_v$.

\begin{coro} \label{coro-main}
Under the hypothesis of the proposition \ref{prop-torsion} on $\mathfrak m$,
the action of $N_\varrho$ on $\Psi_\varrho$ defined above for every 
$\overline \Fm_l$-character $\varrho$, induces a nilpotent 
monodromy operator $N^{coho}_{\mathfrak m,v}$ on
$H^{0}(\sh_{I,\bar s_v},\Psi_{v,\xi})_{\mathfrak m}$ such that
the (multi)-set of Jordan's blocks of  
$\overline N^{coho}_{\mathfrak m,v}:=N^{coho}_{\mathfrak m,v} \otimes_{\overline \Zm_l}
\overline \Fm_l$ acting on $H^{0}(\sh_{I,\bar s_v},\Psi_{v,\xi})_{\mathfrak m} 
\otimes_{\overline \Zm_l} \overline \Fm_l$, is the disjoint union under 
$\{ \widetilde{\mathfrak m} \subset  \mathfrak m \}$, of the
(multi)-sets of Jordan's blocks of $N_{\widetilde{\mathfrak m},v}$.
\end{coro}

\begin{proof}
Recall first that the (multi)-set of Jordan's blocks of
$N^{coho}_{\mathfrak m,v}\otimes_{\overline \Zm_l} \overline \Qm_l$ 
(resp. $\overline N^{coho}_{\mathfrak m,v}$) 
is given by the collection of the dimensions $e_{\mathfrak m,v}(r)$ 
(resp. $\overline e_{\mathfrak m,v}(r)$) of 
$\ker (N^{coho}_{\mathfrak m,v}\otimes_{\overline \Zm_l} \overline \Qm_l)^r$
(resp. $\ker \overline N^{coho}_{\mathfrak m,v}$) for $r \geq 1$: the
columns of the Young diagram associated to 
$N^{coho}_{\mathfrak m,v}\otimes_{\overline \Zm_l} \overline \Qm_l$, 
are of length $e_{\mathfrak m,v}(r+1)-e_{\mathfrak m,v}(r)$ for $r\geq 0$.
\begin{itemize}
\item Proposition \ref{prop-torsion} gives us that the $\overline \Fm_l$-spectral sequence
of nearby cycles degenerates at $E_1$;

\item while from above, we know that
$\overline e_{\mathfrak m,v}(r)$ is the sum over the $\overline \Fm_l$-characters
$\varrho$ of the $\overline \Fm_l$-dimension of
$H^0(\sh_{I,\bar s_v}, \Fil^r_!(\Psi_{\varrho,\xi}))_{\mathfrak m}
\otimes_{\overline \Zm_l} \overline \Fm_l$.

\item But this dimension is also the $\overline \Qm_l$-
dimension of $H^0(\sh_{I,\bar s_v}, \Fil^r_!(\Psi_{\varrho,\xi}))_{\mathfrak m}
\otimes_{\overline \Zm_l} \overline \Qm_l$ 
and so equals to $e_{\mathfrak m,v}(r)$.
\end{itemize}
The result then follows from
$$H^0(\sh_{I,\bar s_v}, \Fil^r_!(\Psi_{\varrho,\xi}))_{\mathfrak m}
\otimes_{\overline \Zm_l} \overline \Qm_l \simeq \bigoplus_{\widetilde{\mathfrak m}
\subset \mathfrak m}
H^0(\sh_{I,\bar s_v}, \Fil^r_!(\Psi_{\varrho,\xi}))_{\widetilde{\mathfrak m}},$$
and the fact that $N_{\widetilde{\mathfrak m},v} = N^{coho}_{\mathfrak m,v} 
\otimes_{\Tm_{\xi,\mathfrak m}^S} \Tm_{\xi,\widetilde{\mathfrak m}}^S$.
\end{proof}

\section{Mazur's principle}
\label{para-main}

\begin{defin} (cf. \cite{scholze-LT} \S 5)   \label{defi-typic}
We say that $\mathfrak m$ is KHT-typic if, as 
a $\Tm_{\xi,\mathfrak m}^S[\gal_{F,S}]$-module,
$$H^{d-1}(\sh_{I,\bar \eta},V_{\xi,\overline \Zm_l})_{\mathfrak m} \simeq 
\sigma_{\mathfrak m} 
\otimes_{\Tm_{\xi,\mathfrak m}^S} \rho_{\mathfrak m},$$
for some $\Tm_{\xi,\mathfrak m}^S$-module $\sigma_{\mathfrak m}$ on which 
$\gal_{F,S}$ acts trivially and 
$$\rho_{\mathfrak m}:\gal_{F,S} \longrightarrow 
GL_d(\Tm_{\xi,\mathfrak m}^S)$$ 
is the stable lattice of
$\bigoplus_{\widetilde{\mathfrak m} \subset \mathfrak m} \rho_{\widetilde{\mathfrak m}}$ introduced in the introduction.
\end{defin}

\rem When $\mathfrak m$ is both KHT-typic and verifies the hypothesis of proposition
\ref{prop-torsion}, then $N_{\mathfrak m,v}^{coho}$ induces a monodromy
operator $N_{\mathfrak m,v}$ on $\rho_{\mathfrak m}$.

As explained in \cite{h-t}, the $\overline \Qm_l$-cohomology of $\sh_{I,\bar \eta}$
can be written as
$$H^{d-1}(\sh_{I,\bar \eta},V_{\xi,\overline \Qm_l})_{\mathfrak m} \simeq 
\bigoplus_{\pi \in \AC_{\xi,I}(\mathfrak m)} (\pi^{\oo})^I \otimes V(\pi^\oo),$$
where 
\begin{itemize}
\item $\AC_{\xi,I}(\mathfrak m)$ is the set of equivalence classes of automorphic
representations of $G(\Am)$ with non trivial $I$-invariants and such that its
modulo $l$ Satake's parameters outside $S$ are prescribed by $\mathfrak m$,

\item and $V(\pi^\oo)$ is a representation of $\gal_{F,S}$. 
\end{itemize}
As $\overline \rho_{\mathfrak m}$ is supposed
to be absolutely irreducible, then as explained in chapter VI of \cite{h-t},
if $V(\pi^\oo)$ is non zero, then $\pi$ is a weak transfer of a $\xi$-cohomological 
automorphic representation $(\Pi,\psi)$ of $GL_d(\Am_F) \times \Am_F^\times$
with $\Pi^\vee \simeq \Pi^c$ where $c$ is the complex conjugation.
Attached to such a $\Pi$ is a global Galois representation 
$\rho_{\Pi,l}:\gal_{F,S} \longrightarrow GL_d(\overline \Qm_l)$ which is irreducible.

\begin{thm} (cf. \cite{nekovar-fayad} theorem 2.20) \\
If $\rho_{\Pi,l}$ is strongly irreducible, meaning it remains irreducible when
it is restricted to any finite index subgroup, then $V(\pi^\oo)$ is a semi-simple
representation of $\gal_{F,S}$.
\end{thm}

\rem The Tate conjecture predicts that $V(\pi^\oo)$ is always semi-simple.

\begin{prop} \label{prop-typic}
We suppose that for all $\pi \in \AC_{\xi,I}(\mathfrak m)$, the Galois representation
$V(\pi^\oo)$ is semi-simple. Then $\mathfrak m$ is KHT-typic.
\end{prop}

\begin{proof}
By proposition 5.4 of \cite{scholze-LT} it suffices to deal with 
$\overline \Qm_l$-coefficients. From \cite{h-t} proposition VII.1.8 and the
semi-simplicity hypothesis, then $V(\pi^\oo) \simeq 
\widetilde R_\xi(\pi)^{\bigoplus n(\pi)}$
where $\widetilde R_\xi(\pi)$ is of dimension $d$. We then write
$$(\pi^\oo)^I \otimes_{\overline \Qm_l} R_\xi(\pi)\simeq 
(\pi^\oo)^I \otimes_{\Tm_{\xi,\mathfrak m,\overline \Qm_l}^S} 
(\Tm_{\xi,\mathfrak m,\overline \Qm_l}^S)^d,$$
and $(\pi^\oo)^I \otimes_{\overline \Qm_l} V(\pi^\oo) \simeq
((\pi^\oo)^I)^{\bigoplus n(\pi)} \otimes_{\Tm_{\xi,\mathfrak m,\overline \Qm_l}^S} 
(\Tm_{\xi,\mathfrak m,\overline \Qm_l}^S)^d$
and finally 
$$H^{d-1}(\sh_{I,\bar \eta},V_{\xi,\overline \Qm_l})_{\mathfrak m} \simeq
M \otimes_{\Tm_{\xi,\mathfrak m,\overline \Qm_l}^S} 
(\Tm_{\xi,\mathfrak m,\overline \Qm_l}^S)^d,$$
with $M \simeq \bigoplus_{\pi \in \AC_{\xi,I}(\mathfrak m)} 
((\pi^{\oo})^I)^{\bigoplus n(\pi)}$. The result then follows from
\cite{h-t} theorem VII.1.9 which insures that $R_\xi(\pi) \simeq
\rho_{\widetilde{\mathfrak m}}$, if $\widetilde{\mathfrak m}$ is the prime ideal
associated to $\pi$,
\end{proof}

\begin{defin} \label{defi-depth}
Let $\overline \rho(\mathfrak m)_\bullet$ be the filtration of 
$\rho_{\mathfrak m} \otimes_{\overline \Zm_l} \overline \Fm_l$ defined by its
iterated socle, that is
\begin{itemize}
\item $\overline \rho(\mathfrak m)_0$ is the socle of 
$\rho_{\mathfrak m} \otimes_{\overline \Zm_l} \overline \Fm_l$

\item and for $i \geq 1$, 
$\overline \rho(\mathfrak m)_i/\overline \rho(\mathfrak m)_{i-1}$ is
the socle of $(\rho_{\mathfrak m} \otimes_{\overline \Zm_l} \overline \Fm_l)
/\overline \rho(\mathfrak m)_{i-1}$.
\end{itemize}
The depth of $\mathfrak m$ is then the length of this filtration.
\end{defin}

\begin{thm} \label{theo-main} (\textbf{Mazur's principle})
Let $\mathfrak m$ be a maximal ideal of $\Tm_{\xi,\overline \Zm_l}^S$ such
that:
\begin{itemize}
\item $\overline \rho_{\mathfrak m}$ is absolutely irreducible and its restriction
to the decomposition group at $v$ is, up to the action of the monodromy operator,
a direct sum of characters;

\item $\mathfrak m$ is KHT-free and KHT-typic.

%\item $H^{d-1}(\sh_{I,\bar \eta},V_{\xi,\overline \Qm_l})_{\mathfrak m}$ is Galois semi-simple.
\end{itemize}
Let $\underline{\bar d_{\mathfrak m,v}}=(t_1 \geq \cdots \geq t_r)$ 
be the partition of $d$ given by the Jordan blocks of $\overline N_{\mathfrak m,v}$.
Then there exists $\widetilde{\mathfrak m} \subset \mathfrak m$ such that
$$\rho_{\widetilde{\mathfrak m},v} \simeq \sp_{t_1}(\chi_{v,1}) \oplus \cdots \oplus
\sp_{t_r}(\chi_{v,r}),$$
where $\chi_{v,i}$ are non isomorphic characters. 
\end{thm}

For $\underline d$ a partition of $d$, its associated parahoric subgroup is
\addtocounter{thm}{1}
\begin{equation} \label{eq-parahoric}
I_{\underline d}(\OC_v)=\ker (GL_d(\OC_v) \longrightarrow P_{\underline d}(\kappa(v)))
\end{equation}
where $P_{\underline d}$ is the standard parabolic subgroup associated
to $\underline d$. The dominance order $(d_1 \geq \cdots d_r) \leq (e_1 \geq \cdots \geq
e_s)$ is the given by 
$$\forall i \geq 1: ~\sum_{i=1}^k d_i \leq \sum_{i=1}^k e_i.$$
Recall, cf. lemma 1.1.7 of \cite{boyer-mazur} that $\Pi_v\simeq \st_{t_1}(\chi_{v,1}) \times \cdots \times \st_{t_r}(\chi_{v,r})$
has non trivial invariant vectors upon $I_{\underline d}(\OC_v)$ if and only if
$\underline d$ is smaller than the dual partition $\underline e^*$ of 
$\underline e:=(t_1 \geq \cdots \geq t_r)$ whose lines are the rows of $\underline e$.
The theorem then says that $\Pi_{\widetilde{\mathfrak m},v}$ has non trivial
invariant vectors upon $I_{\underline{\bar d_{\mathfrak m,v}}^*}(\OC_v)$.

\begin{proof}
We will consider three different types of situations at $v$:
\begin{itemize}
\item infinite and we then denote by $H(I^v(\oo))_{\mathfrak m}:=
H^0(\sh_{I^v(\oo),\bar s_v},\Psi_{v,\xi})_{\mathfrak m}$ and
$\overline H(I^v(\oo))_{\mathfrak m}:= H(I^v(\oo))_{\mathfrak m}
\otimes_{\overline \Zm_l} \overline \Fm_l$;

\item $H(I^v,\underline d)_{\mathfrak m}:=
H^0(\sh_{I^vI_{\underline d}(\OC_v),\bar s_v},\Psi_{v,\xi})_{\mathfrak m}$
and $\overline H(I^v,\underline d)_{\mathfrak m}:=H(I^v,\underline d)_{\mathfrak m}
\otimes_{\overline \Zm_l} \overline \Fm_l$;

\item $H(I^v,h)_{\mathfrak m}:=J^{GL_d}_{P_h} \bigl ( H(I^v(\oo))_{\mathfrak m} 
\bigr )^{GL_h(\OC_v)}$ and $\overline H(I^v,h)_{\mathfrak m}:=
H(I^v,h)_{\mathfrak m} \otimes_{\overline \Zm_l} \overline \Fm_l$; where
$J^{GL_d}_{P_h}$ is the Jacquet functor associated to the parabolic
$P_h(F_v)$.
\end{itemize}
We also denote by $\Tm_{\xi,\mathfrak m}(I^v,\underline d)$ the image of 
$\Tm^S_{abs}$ inside $H(I^v,\underline d)_{\mathfrak m}$.
Let $\underline d$ be minimal for the previous dominance order, such that 
$H(I^v,\underline d^*)_{\mathfrak m} \neq 0$ and
consider now a broken row in $\overline T_{\mathfrak m,v}$ that is
$\overline \lambda \in \overline \Fm_l$ such that $q_v \overline \lambda$
does not belong to the same line with $\overline \lambda$ 
in the labelled Young diagram of $\overline N_{\mathfrak m,v}$.

\begin{lemma} \label{lem-coupure}
There exists $\widetilde{\mathfrak m} \subset \mathfrak m$ with
$H(I^v,\underline d^*)_{\widetilde{\mathfrak m}} \neq (0)$, such that 
inside $T_{\widetilde{\mathfrak m},v}$, 
the liftings\footnote{cf. the multiplicity free hypothesis} $\lambda_1$ and $\lambda_2$
of $q_v \overline \lambda$ and $\overline \lambda$, do not belong to the same line, 
i.e. $\lambda_1/\lambda_2 \neq q_v^{\pm 1}$.
\end{lemma}

\begin{proof}
Let denote by $\rho_{\mathfrak m}(\underline d):=\rho_{\mathfrak m} \otimes_{\Tm_\xi^S}
\Tm_\xi(I^v,\underline d^*)$ and $N_{\mathfrak m,v}(\underline d):=N_{\mathfrak m,v}
\otimes_{\Tm_\xi^S} \Tm_\xi(I^v,\underline d^*)$.
We consider then the eigenspaces $V_I(\overline \lambda)$ and
$V(\overline \lambda)$ for the eigenvalue 
$\overline \lambda$ of the action of $\frob_v$, respectively, 
on the $\overline \Fm_l$-vector spaces
$$\im N_{\mathfrak m,v}(\underline d) \otimes_{\overline \Zm_l} \overline \Fm_l \subset \rho_{\mathfrak m}(\underline d) \otimes_{\overline \Zm_l} \overline \Fm_l.$$ 
Note that, as the eigenvalues of $\frob_v$ acting on
$\rho_{\mathfrak m} \otimes_{\overline \Zm_l} \overline \Fm_l$ are supposed
to be pairwise distinct, then 
$$\dim_{\overline \Fm_l} V(\overline \lambda)=
\sharp \{ \widetilde{\mathfrak m} \subset \mathfrak m; , H(I^v,\underline d^*)_{\widetilde{\mathfrak m}} \neq (0) \}.$$

Note that if the conclusion of the lemma were not true, then
$\dim_{\overline \Fm_l} V_I(\overline \lambda)=\dim_{\overline \Fm_l} 
V(\overline \lambda)$. Indeed this is true over
$\overline \Qm_l$ and from corollary \ref{coro-main},  
we know that the image of 
$N_{\mathfrak m,v}^{coho} \otimes_{\overline \Zm_l} \overline \Fm_l$
is the modulo $l$ reduction of the image of $N_{\mathfrak m,v}$, which gives us
the previous equality.

Consider now the previous filtration $\overline \rho(\mathfrak m,\underline d)_i$ of
$\rho_{\mathfrak m}(\underline d) \otimes_{\overline \Zm_l} \overline \Fm_l$, where all the
graded parts $\gr_i(\mathfrak m,\underline d)$ are a direct sum of $\overline \rho_{\mathfrak m}$.
If we want $\dim_{\overline \Fm_l} V_I(\overline \lambda)$ to be equal to the number 
$\sharp \{ \widetilde{\mathfrak m} \subset \mathfrak m; , H(I^v,\underline d^*)_{\widetilde{\mathfrak m}} \neq (0) \}$ of irreducible subquotients
of $\rho_{\mathfrak m}(\underline d) \otimes_{\overline \Zm_l} \overline \Fm_l$, 
then $\frob_v$  should induces an isomorphism 
$V(q_v \overline \lambda) \rightarrow V_I(\overline \lambda)$.
But note that $V(q_v \overline \lambda)$ intersects 
$\overline \rho(\mathfrak m,\underline d)_0$
and, as $q_v \overline \lambda \rightarrow \overline \lambda$ is broken in 
$\overline T_{\mathfrak m,v}$, and the eigenvalues are all distinct, 
then $\frob_v$ acts trivially on
$V(q_v \overline \lambda) \cap \overline \rho(\mathfrak m)^0$ so that
$\dim_{\overline \Fm_l} V_I(\overline \lambda) < \dim_{\overline \Fm_l} 
V(q_v \overline \lambda)=
\sharp \{ \widetilde{\mathfrak m} \subset \mathfrak m;, H(I^v,\underline d^*)_{\widetilde{\mathfrak m}} \neq (0) \}$.
\end{proof}

\noindent \textit{End of the proof}: recall that $\underline d$ was taken 
minimal such that  $H(I^v,\underline d^*)_{\mathfrak m} \neq 0$. 
By freeness of the cohomology groups, there exists $\widetilde{\mathfrak m}
\subset \mathfrak m$ with 
$H_{\overline \Qm_l}(I^v,\underline d^*)_{\mathfrak m} \neq 0$, so that
$\underline d_{\widetilde{\mathfrak m},v}$ is greater that $\underline d$ and finally,
by minimality of $\underline d$, we have 
$\underline d=\underline d_{\widetilde{\mathfrak m},v}$, and it remains to prove
that they are equal to $\underline d_{\mathfrak m,v}$.
Note that $\underline d_{\widetilde{\mathfrak m},v}$ is obtained from
$\underline d_{\mathfrak m,v}$ by glueing some of its line so that if it were note
equal to $\underline d_{\mathfrak m,v}$, the previous lemma would tell us that
$\underline d$ were not maximal.

%
% we then consider
%$\overline H(I_v,h)_{\mathfrak m}$ where $h$ is the length of the first
%column in the Young diagram associated to $\underline d$. Note that
%$\overline H(I_v,h)_{\mathfrak m}$ is also a $\Tm_{h,v}$-module where
%$\Tm_{h,v}$ is the unramified Hecke algebra of $GL_h(\OC_v)$.
%For a maximal ideal $\mathfrak m_v(h)$ of $\Tm_{h,v}$, we can
%consider $\overline H(I_v,h)_{\mathfrak m(h,v)}$ where 
%$\mathfrak m(h,v):=\mathfrak m \otimes \mathfrak m_v(h)$ in
%$\Tm_{\xi,\mathfrak m}^S \otimes \Tm_{h,v}$.
%
%Note that $\widetilde{\mathfrak m}$ contributes to 
%$H(I^v,h)_{\mathfrak m}$ (resp. $H(I^v,h)_{\mathfrak m(h,v)}$)
%if and only if the first column of the Young 
%diagram associated to $N_{\widetilde{\mathfrak m},v}$ is of length greater than $h$
%(resp. the labels in the first columns contained the prescribed sets given by
%$\mathfrak m_v(h)$). As $\underline d$ is taken maximal,

\end{proof}

%The following result was explained to me by one of the anonymous referee.

\begin{prop}
With the previous notations, for every $0 \leq i \leq r-1$ where $r$ is the depth of 
$\mathfrak m$, then $\gr^i(\mathfrak m)$ is irreducible and then isomorphic to
$\overline \rho_{\mathfrak m}$.
\end{prop}

\begin{proof}
Automorphic representation $\pi_{\widetilde{\mathfrak m}}$ in level $I^v(\oo)$ having
a cuspidal support at $v$ made of characters, are in finite numbers so that
$\rho_{\mathfrak m} \otimes_{\overline \Zm_l} \Tm_{\mathfrak m,\xi}$
is defined over a noetherian ring $R$ inside 
$\prod_{\widetilde{\mathfrak m} \subset \mathfrak m} K_{\widetilde{\mathfrak m}}$
where the finite extension $K_{\widetilde{\mathfrak m}}/\Qm_l$ is the field of definition of
$\rho_{\widetilde{\mathfrak m}}$. We then
denote by $\kappa=R/\mathfrak m$ and we still denote by 
$\overline \rho_{\mathfrak m}$ the associated $\kappa$-vector space with its
action of the Galois group.

By construction $\gr^i(\mathfrak m)$ is then a semi-simple $\gal_{F,S}$-module
with underlying representation space a free rank
$\Tm_{\xi,\mathfrak m}^{S,(i)} \otimes_R \kappa$-module $V^{(i)}$ where
$\Tm_{\xi,\mathfrak m}^S \twoheadrightarrow \Tm_{\xi,\mathfrak m}^{S,(i)}$.
We then have a $\Tm_{\mathfrak m,\xi}^{S,(i)}[\gal_{F,S}]$-equivariant isomorphism
$$\overline \rho_{\mathfrak m} \otimes_R \kappa \simeq
\overline \rho_{\mathfrak m} \otimes_{\kappa} \hom_{\gal_{F,S}}
(\overline \rho_{\mathfrak m},V),$$
where $\gal_{F,S}$ acts on the first factor, and $\Tm_{\mathfrak m,\xi}^{S,(i)}$
on the second. As a deformation of $\overline \rho_{\mathfrak m}$ over 
the local $\kappa$-algebra $\Tm_{\xi,\mathfrak m}^{S,(i)} \otimes_R \kappa$,
it is given by a map
$R_{\overline \rho_{\mathfrak m}} \longrightarrow 
\Tm_{\xi,\mathfrak m}^{S,(i)} \otimes_R \kappa$ where 
$R_{\overline \rho_{\mathfrak m}}$ is the usual deformation ring. The previous
isomorphism implies that this last map has to factor through the residue field
of $R_{\overline \rho_{\mathfrak m}}$ which is $\kappa$.
Moreover from $R=T$ theorem, we have a surjection 
$$R_{\overline \rho_{\mathfrak m}} \twoheadrightarrow T_{\xi,\mathfrak m}^S
\twoheadrightarrow T_{\xi,\mathfrak m}^{S,(i)}$$ 
so that
$T_{\xi,\mathfrak m}^{S,(i)} \simeq \kappa$. By Nakayama's theorem, 
$T_{\xi,\mathfrak m}^{S,(i)}$ is then a local ring corresponding to an unique
$\widetilde{\mathfrak m}$.
\end{proof}

The set of partitions $\underline d_{\widetilde{\mathfrak m},v}$
for various $\widetilde{\mathfrak m} \subset \mathfrak m$, could be used
to obtain informations about the depth of $\mathfrak m$. Consider for example the following situation:
\begin{itemize}
\item $S_v(\mathfrak m)= \{ \alpha,q_v \alpha,\cdots, q_v^{d-1} \alpha \}$;

\item $\overline N_{\mathfrak m,v}$ is zero;

\item there exists $\widetilde{\mathfrak m} \subset \mathfrak m$ such that
$\Pi_{\widetilde{\mathfrak m},v} \simeq \st_d(\chi_v)$ for some character
$\chi_v$.
\end{itemize}

\begin{lemma}
With the hypothesis of \ref{theo-main} and the three above assumptions, then
the depth of $\mathfrak m$ is greater than $d$.
\end{lemma}

\begin{proof}
By construction each of the $\overline \rho_i(\mathfrak m)$ is a direct sum of
copies of $\overline \rho_{\mathfrak m}$ so that the nilpotent monodromy operator
$\overline N_{\mathfrak m,v}$ acts trivially. We then deduce that 
$\overline N_{\mathfrak m,v}$ sends $\overline \rho(\mathfrak m)_i$ onto
$\overline \rho(\mathfrak m)_{i-1}$. Our last hypothesis then implies that
$\overline N_{\mathfrak m,v}^{d-1} \neq 0$ so that the depth of
$\mathfrak m$ should be greater than $d$.
\end{proof}

%As explained in \S \ref{para-ihara}, the existence of 
%$\widetilde{\mathfrak m} \subset \mathfrak m$ such that
%$\Pi_{\widetilde{\mathfrak m},v} \simeq \st_d(\chi_v)$,
%should be given by the
%higher dimension of Ihara's lemma. 
More generally, consider 
\begin{itemize}
\item $r$ maximal such that there exists $\alpha$ with
$\{ \alpha,q_v \alpha, \cdots,q_v^{r-1} \alpha \} \subset S_v(\mathfrak m)$.
We also denote by $e_0,\cdots,e_{r-1}$ the associated eigenvectors
of $\overline \rho_{\mathfrak m}(\frob_v)$.

\item Let denote by $i_0=0 < i_1 < \cdots < i_k \leq r-1$ the indexes $i$ such that
$e_i \in \ker \overline N_{\mathfrak m,v}$.

\item We moreover assume the existence of $\widetilde{\mathfrak m} \subset
\mathfrak m$ such that $\Pi_{\widetilde{\mathfrak m}} \simeq \st_r(\chi_v) \times ?$
where $\chi_v$ is a character of $F_v^\times$ such that $\chi_v(\varpi_v)=\alpha$
and where $?$ means a irreducible representation we do not want to precise.
\end{itemize}

\begin{lemma}
With the hypothesis of \ref{theo-main} and the three above assumptions, then
the depth of $\mathfrak m$ is strictly greater than $k$.
\end{lemma}

\begin{proof}
The existence of $\widetilde{\mathfrak m}$ implies that there exists
an eigenvector $f_{r-1}$ of $\rho_{\mathfrak m}(\frob_v) 
\otimes_{\overline \Zm_l} \overline \Fm_l$ for the eigenvalue $q_v^{r-1} \alpha$
such that $(\overline N_{\mathfrak m,v}^{coho})^{r-1}(f_{r-1}) \neq 0$.
We first introduce the following notations:
\begin{itemize}
\item $i$ such that $f_{r-1} \in \overline \rho(\mathfrak m)_i$;

\item for $ \leq j \leq r-1$, let $f_j=N_{\mathfrak m,v}^{r-1-j}(f_{r-1})$.
\end{itemize}
Through $\overline \rho(\mathfrak m)_i \twoheadrightarrow \overline \rho(\mathfrak m)_i
/  \overline \rho(\mathfrak m)_{i-1}$, 
the image of $f_{i_k}$ by $\overline N_{\mathfrak m,v}$ is zero so that 
$f_{i_k-1} \in \overline \rho(\mathfrak m)_{i-1}$.
As $S_v(\mathfrak m)$ is supposed to be multiplicity free then the
image of $f_{i_k-1}$ in $\overline \rho(\mathfrak m)_{i-1}/\overline \rho(\mathfrak m)_{i-2} 
\simeq \overline \rho_{\mathfrak m}^{\bigoplus m_{i-1}}$ belongs to the space
generated by the $e_{i_k-1}$ in each of the copies of 
$\overline \rho_{\mathfrak m}$. We can then repeat the previous observation so
that the image of $f_{i_{k-1}-1} \in \overline \rho(\mathfrak m)_{i-2}$ 
and that finally the depth of $\mathfrak m$ should be
greater than $k$.
\end{proof}

\section{The reducible case}
\label{para-reducible}

In \cite{boyer-mazur}, $N_{\mathfrak m,v}^{coho}
\otimes_{\overline \Zm_l} \overline \Fm_l$
was wrongly identified with the direct sum of $\overline N_{\mathfrak m,v}$, so that
\cite{boyer-mazur} proposition 3.1.12 is false. The aim of loc. cit. was to give conditions 
on $\underline{d_{\mathfrak m,v}}$ so that the torsion of the cohomology of the
KHT-Shimura variety is non trivial, and then play with the level in order to obtain
level lowering. 

In this section we want to resume the strategy of \cite{boyer-mazur} 
and give a correct statement about level lowering.
From the main result of \cite{boyer-jep2}, we then have to consider
the case where $\overline \rho_{\mathfrak m}$ is reducible. In \cite{boyer-mazur} 
the irreducibility hypothesis for $\overline \rho_{\mathfrak m}$ was used to insure
the free quotients of the various cohomology groups to be concentrated in
middle degree. To keep this property we then make the following hypothesis.

\medskip

\noindent \textit{Hypothesis}:
We now consider the case where 
\begin{itemize}
\item $\overline \rho_{\mathfrak m}$ is reducible 

\item and the semi-simplification of its restriction to the decomposition group at $v$, 
which we write 
$\bigoplus_{i=1}^r \sp_{s_i}(\chi_{v,i})$ without taking into account the monodromy
operator, verifies the following property where 
$\sp_{s_i}(\chi_{v_i})=\chi_{v,i} (\frac{1-s_i}{2}) \oplus \cdots  \chi_{v,i} (\frac{s_i-1}{2})$ 
as before:
\begin{itemize}
\item the $\chi_{v,i}(\frac{1-s_i+2t_i}{2})$ for $i=1,\cdots,r$ and 
$0 \leq t_i \leq s_i-1$, are pairwise distinct, i.e.
$\chi_{v,i}(\frac{1-s_i+2t_i}{2})$ is not isomorphic to any of the
$\chi_{v,j}(\frac{1-s_j+2t_j}{2})$ whatever are $t_i$ and $t_j$ verifying the previous
inequality. Taking $i=j$, we obtain in particular that,
for all $i=1,\cdots,r$, the order of $q_v$ modulo $l$ is strictly greater than $s_i$.

\item Finally we impose that there exists $i \neq j$ with $s_i \neq s_j$.
\end{itemize}
\end{itemize}

Recall that if $\widetilde{\mathfrak m} \subset \mathfrak m$ appears in the cohomology
of our KHT Shimura variety such that $\rho_{\widetilde{\mathfrak m}}$ is 
reducible then $T_{\widetilde{\mathfrak m},v}$ is a rectangle. It is in particular
excluded from the previous hypothesis, i.e. 
for every $\widetilde{\mathfrak m} \subset \mathfrak m$
then $\rho_{\widetilde{\mathfrak m}}$ is irreducible. 
In particular the free quotients
of the cohomology localized at $\mathfrak m$ are still concentrated in middle
degree but there might be non trivial torsion classes.
Note also that the modulo $l$ reduction of $\rho_{\widetilde{\mathfrak m}}$
depends on $\widetilde{\mathfrak m}$ but also of the chosen stable lattice
inside it.

As now KHT-typicness is no more verified, we consider, for a partition $\underline d$
of $d$, a $\gal_{F,S}$-equivariant filtration of the free quotient 
$H(I^v,\underline d^*)_{\mathfrak m,free}$ of $H(I^v,\underline d^*)_{\mathfrak m}$,
so that each graded part is a stable lattice of some $\rho_{\widetilde{\mathfrak m}}$
with $\pi_{\widetilde{\mathfrak m},v}^{I_{\underline d^*}(\OC_v)} \neq (0)$.
The modulo $l$ reductions of each of these lattices, give then partitions
of $d$ associated to the Jordan blocks of the monodromy operator at $v$ and
we denote by $\underline{d_{\mathfrak m,v}}(\underline d)$ the minimal one.

\rem Note that, a priori, $\underline{d_{\mathfrak m,v}}(\underline d)$ 
might depend on the choice of the filtration of 
$H(I^v,\underline d^*)_{\mathfrak m,free}$.

\begin{defin} (cf. \cite{boyer-mazur} \S 1.1 and definition 1.3.1) \\
For $\underline d$ a partition of $d$, we denote by $\underline d^{(1)}$
the partition such that its Young diagram is obtain from those of
$\underline d$ by deleting its first column.

We then say that $\mathfrak m$ is degraded relatively to $\underline d$ if
$\underline{d}^{(1)}=(d_1 \geq d_2 \geq \cdots )$ is not contained in
$\underline{d_{\mathfrak m,v}}(\underline d)=(t_1 \geq t_2 \geq \cdots)$ i.e.
there exists $i \geq 1$ such that $d_i \geq t_i$.
\end{defin}

\begin{prop}
Let $\underline d$ minimal such that 
$H(I^v,\underline d^*)_{\mathfrak m,free} \neq (0)$ and suppose that
$\mathfrak m$ is degraded relatively to $\underline d$, then 
for every $w \in \spl(I) \setminus \{ v \}$, there exists
$\widetilde{\mathfrak m} \subset \mathfrak m$ with non zero invariant vectors upon
$I(w)$ defined as follows:
\begin{itemize}
\item outside $v$ and $w$ it coincides with $I^{v,w}$;

\item at $w$ it is of parahoric type for a partition $(t,1,\cdots, 1)$ for some $t \leq d$;

\item at $v$, $I(w)_v=I_v(\underline e^*)$ with $\underline e < \underline d$.
\end{itemize}
\end{prop}

\begin{proof}
By hypothesis
$\underline{d_{\mathfrak m,v}}(\underline d)$ is strictly smaller than $\underline d$.
We first want to prove that in the spectral sequence which degenerates in $E_1$ over
$\overline \Qm_l$, 
$$E_{1,\mathfrak m}^{p,q}=H^{p+q}(\sh_{I^vI_{\underline d^*}(\OC_v),\bar s_v},
\grr^{-p}(\Psi_{v,\xi}))_{\mathfrak m} \Rightarrow 
H^{p+q}(\sh_{I^vI_{\underline d^*}(\OC_v),\bar s_v},\Psi_{v,\xi})_{\mathfrak m},$$
some of the $E_{1,\mathfrak m}^{p,q}$ for $p+q \in \{ 0,1 \}$ 
are not torsion free. We argue by contradiction. From lemma \ref{lem-ih},  there
is then no nontrivial torsion classes in the initial terms of the spectral sequence.
We follows the proof of theorem \ref{theo-main}. First the conclusion of 
lemma \ref{lem-coupure} is still true. Indeed 
\begin{itemize}
\item consider a $\gal_{F,S}$-equivariant filtration $\Fil^\bullet$ of 
$H(I^v,\underline d^*)_{\mathfrak m,free}$ 
such that the graded parts $\gr^k$ are free and irreducible after inverting $l$.

\item Modulo $l$ we then obtain a filtration with graded parts $\overline \gr^k$ isomorphic 
to $\Gamma/l\Gamma$ where $\Gamma$ is some stable lattice of some 
$\rho_{\widetilde{\mathfrak m}}$: we then denote by $\underline d(k)$ the associated
partition given by the monodromy operator.

\item Let $k$ be minimal such that in the labelled Young diagram of $\underline d(k)$,
there exists
$q_v \bar \lambda$ and $\bar \lambda$  which are not in the same line: 
by hypothesis such a $k$ exists. 

\item Let denote by $V_k(\bar \lambda)$ the eigenspace
of $\frob_v$ acting on $\overline \Fil_k$ for the eigenvalue $\bar \lambda$.
Note then that the dimension of 
$\overline N^{coho}_{\mathfrak m,v} V_k(q_v \bar \lambda)$ is strictly less than 
$\dim V_k(\bar \lambda)$ and this inequality remains true for every $k' \geq k$.

\item The conclusion is then similar to those in the proof of  \ref{lem-coupure}.
\end{itemize}
If there were no torsion, then following the proof in the previous section, 
by minimality of $\underline d$, we would obtain
$\widetilde{\mathfrak m} \subset \mathfrak m$ such that 
$\underline d=\underline{d_{\widetilde{\mathfrak m},v}}$ would be equal to 
$\underline{d_{\mathfrak m,v}}(\underline d)$
which is not the case as  $\mathfrak m$ is supposed to be degraded relatively to
$\underline d$.

So we know that at least one of the $E_{1,\mathfrak m}^{p,q}$ for $p+q \in \{ 0,1 \}$, 
has non trivial torsion classes. 
%
%
%The idea is then to follow the proof of proposition \ref{prop-torsion}
%to prove that $E_{\oo,\mathfrak m}^i$ has some non trivial torsion classes.
% About the hypothesis of \ref{prop-torsion}
% \begin{itemize}
% \item note first that the hypothesis on the irreducibility of $\overline \rho_{\mathfrak m}$
%in the proof of \ref{prop-torsion} is only use to insure that the free quotients on the
%cohomology groups are concentrated in middle degree. This property remains
%true here as already noticed before.
%
%\item The third hypothesis is also verified from our assumptions.
%\end{itemize}
%By contraposition, proposition \ref{prop-torsion} then gives us the  existence of
%non trivial torsion classes in the cohomology groups $H^i(\sh_{I,\bar \eta_v},
%V_{\xi,\overline \Zm_l})_{\mathfrak m}$. 
We now take up the arguments of \cite{boyer-mazur} \S 3.3:
\begin{itemize}
\item There should exists $h_0$ such that the cokernel of (\ref{eq-map1-coho2}) has 
non trivial torsion. In particular, cf. \cite{boyer-compositio}, there should exists 
an automorphic representation $\Pi$, irreducible, $\xi$-cohomological with
$$\Pi_v \simeq \st_{t_0+1}(\chi_{v,0}) \times \st_{t_1}(\chi_{v,1}) \times \cdots \times
\st_{t_r}(\chi_{v,r}),$$
where $\chi_{v,i}$ are characters of $F_v^\times$ and $t_1,\cdots,t_r$
are integers we do not want to precise at this point;

\item to have non trivial torsion in level $I^vI_v(\underline d^*)$, $\Pi^v$ should have
non trivial invariant vectors under $I^v$ and the partition associated to
$(t_0,1,t_1,\cdots, t_r)$ should be less than $\underline d$. 

\item From lemma \ref{lem-ih}, then $H^{1-t_0}(\sh_{I^vI_v(\underline e^*),\bar s_v},
\lexp p j^{=1}_{!*} HT(\chi_v,\chi_v))_{\mathfrak m}$ has non trivial torsion subspace
where $\underline e^*$ is the dual of 
$(\overbrace{1,\cdots,1}^{t_0+1},t_1, \cdots,t_r)$.

\item If $t_0 \geq 2$ then as $(\overbrace{1,\cdots,1}^{t_0+1},t_1, \cdots,t_r) <
(t_0,1,t_1,\cdots,t_r) \leq \underline d$. We then consider the spectral sequence
(\ref{eq-ss-gr}) in level $I^vI_v(\underline e^*)$ so that the initial terms are all torsion, some
of then being non zero. The proof of proposition \ref{prop-torsion} then gives us
that $E_\oo^{i}$ are then torsion and non zero for some $i$.
\end{itemize}

We now have to deal with the case where $t_0=1$ and $\underline e=\underline d$.
As in the proof of 
\cite{boyer-mazur} lemma 3.3.3, we then consider
$\coFil^1_!(\Psi_{v,\xi}):=\Psi_{v,\xi}/\Fil^1_!(\Psi_{v,\xi})$ and its 
$\overline \Zm_l$-exhaustive filtration where now $\PC(1,\chi_v)$ does not appear
anymore so that the spectral sequence associated to this filtration degenerates
as before in $E_1$ with torsion free initial terms. In particular as in the proof
of corollary \ref{coro-main}, the dimensions $e_{\mathfrak m,v}(r)$ (resp. $\overline
e_{\mathfrak m,v}(r)$)
of $N_{\mathfrak m,v}^r \otimes_{\overline \Zm_l} \overline \Qm_l$ (resp.
$(\overline N_{\mathfrak m,v}^{coho})^r$) for $r \geq 2$ coincide while for $r=1$
we have $\overline e_{\mathfrak m,v}(1)=e_{\mathfrak m,v}(1)+\delta$ for some
$0 \leq \delta \leq e_{\mathfrak m,v}(2)-e_{\mathfrak m,v}(1)$. 

The labelled Young's diagrams of $\overline N_{\mathfrak m,v}^{coho}$ are 
then obtained from those of $\underline d$ allowing to untie its first column. 
We then argue as before to conclude that $\underline d^{(1)}$ has to be contained
in every labelled Young's diagram of the modulo $l$ reduction, relatively to
lattices given by the cohomology, of $\rho_{\widetilde{\mathfrak m,v}}$
and so in particular $\underline d^{(1)}$ is contained in
$\underline{d_{\mathfrak m,v}}(\underline d)$ telling that $\mathfrak m$
is not degraded.

Finally, we deduce that there exists some non trivial torsion classes in some of the
$H^i(\sh_{I^vI_v(\underline d'),\bar \eta_v}, V_{\xi,\overline \Zm_l})_{\mathfrak m}$. 
From the main result of \cite{boyer-mrl}, up to increase the level at some extra place
$w \in \spl(I) \setminus \{ v \}$ where the level become parahoric for
$(t,1\cdots,1)$ for some $t \leq d$, we can lift $\mathfrak m$ with level
$I_v(\underline e^*)$ at $v$.

\end{proof}

\section{Automorphic congruences}
\label{para-congruence}

As in  \cite{boyer-aif}, we can use the freeness of the cohomology groups
of the Harris-Taylor perverse sheaves, to produce automorphic congruences.
Consider then $\mathfrak m$ verifying the hypothesis of proposition 
\ref{prop-torsion} so that the 
$H^i(\sh_{I,\bar s_v},\lexp p j^{=h}_{!*}HT_\xi(\chi_v,h))_{\mathfrak m}$ are free 
and concentrated in degree $i=0$ with
\addtocounter{thm}{1}
\begin{multline} \label{eq-equality}
H^0(\sh_{I,\bar s_v},\lexp p j^{=h}_{!*}HT_\xi(\chi_v,h))_{\mathfrak m} 
\otimes_{\overline \Zm_l} \overline \Fm_l \simeq \\
H^0(\sh_{I,\bar s_v},\lexp p j^{=h}_{!*} HT_{\xi,\overline \Fm_l}(r_l(\chi_v),h) )_{\mathfrak m} \\ 
\simeq H^0(\sh_{I,\bar s_v},\lexp p j^{=h}_{!*} HT_\xi(\chi'_v,h))_{\mathfrak m} 
\otimes_{\overline \Zm_l} \overline \Fm_l,
\end{multline}
whatever is $\chi'_v$ such that the modulo $l$ reduction $r_l(\chi'_v)$ of $\chi'_v$
is isomorphic to those of $\chi_v$. Recall then from \cite{boyer-compositio},
the description of the $\overline \Qm_l$-cohomology groups of 
$\lexp p j^{=h}_{!*} HT_\xi(\chi_v,h)$ localized at $\mathfrak m$.

\begin{prop} \label{prop-cohomo}
(cf. \cite{boyer-compositio} \S 3.6 with\footnote{As $\overline \rho_{\mathfrak m}$ is supposed to be irreducible, the integer
$s$ of \cite{boyer-compositio} \S 3.6 is necessary equal to $1$.}  $s=1$) \\
For $\chi_v$ an unitary character of $F_v^\times$, then, for
$1 \leq h \leq d$, as a $\Tm_{\mathfrak m}^S[GL_d(F_v)]$-module, we have
$$\lim_{\atop{\rightarrow}{I_v}} 
H^0(\sh_{I^vI_v,\bar s_v},\lexp p j^{=h}_{!*} HT_{\xi,\overline \Qm_l}
(\chi_v,h))_{\mathfrak m} \simeq
\bigoplus_{\Pi \in \AC(I,\xi,h,\chi_v,\mathfrak m)} m(\Pi) (\Pi^{\oo,v})^{I^v} \otimes \Pi_v,$$
where 
\begin{itemize}
\item $\AC(I,h,\chi_v,\mathfrak m)$ is the set of irreducible $\xi$-cohomological
automorphic representations $\Pi$ of $G(\Am)$ with non zero invariants under $I^v$
with modulo $l$ Satake's parameters prescribed by $\mathfrak m$,

\item such that $\Pi_v$ is of the following shape 
$$\Pi_v\simeq \st_h(\chi_v) \times \Psi_v$$ 
where $\Psi_v$ is a representation of $GL_{d-h}(F_v)$,

\item and $m(\Pi)$ is the multiplicity of $\Pi$ in the space of automorphic forms.
\end{itemize}
\end{prop}

\rem We write the local component $\Pi_v$ of 
$\Pi \in  \AC(I,\xi,h,\chi_v,\mathfrak m)$ as 
$$\Pi_v \simeq \st_{t_1}(\chi_{v,1}) \times \cdots \times \st_{t_r}(\chi_{v,r}) \times \Psi'_v,$$
where 
\begin{itemize}
\item the $\chi_{v,i}$ are inertially equivalent characters, 
%i.e. $\chi_{v,i} \simeq \chi_{v,1}\otimes \xi_i \circ \val$ for characters $\xi_i:\Zm \longrightarrow \overline \Zm_l^\times$: we then write $\chi_{v,i} \sim_i \chi_{v,1}$;

\item $\Psi'_v$ is an irreducible representation of $GL_{d-\sum_{i=1}^r t_i}(F_v)$ 
whose cuspidal support, made of character by hypothesis, 
does not contain a character inertially equivalent to $\chi_{v,1}$.
\end{itemize}
Then $\Pi$ contributes $k$ times in the isomorphism of the previous proposition, where
$k=\sharp \{ 1 \leq i \leq r \hbox{ such that } t_i=h \}$.

%
%\begin{proof}
%From \cite{boyer-torsion}, there exists a filtration of stratification 
%$$(0)=\Fil^0(\chi_{v} ,h) \hookrightarrow \Fil^{-d}(\chi_{v} ,h) \hookrightarrow \cdots \\
%\hookrightarrow \Fil^{-h}(\chi_{v} ,h)=j^{=h}_{!} HT(\chi_{v} ,\Pi_h)$$
%such that, using lemma \ref{lem-ext0}, the graded parts are 
%$$\gr^{-k}(\chi_{v},h) \simeq \lexp p j^{=k}_{!*} 
%HT(\chi_{v},\Pi_h \{\frac{h-k}{2} \} \otimes \st_{k-h}(\chi_{v}\{h/2 \} ))(\frac{h-k}{2}).$$ 
%The $\xi$-associated spectral sequence localized at $\mathfrak m^{v}$ is then concentrated in 
%middle degree and torsion free.
%Then the spectral sequence associated to this filtration has all its $E_1$ 
%terms torsion free and degenerates at its $E_1$ terms. 
%\end{proof}

%Note also
%that $\Psi_v$ in the statement is of the following shape
%$\st_{t_1}(\pi_{v,1}) \times \cdots \times \st_{t_r}(\pi_{v,r})$ where the $\pi_{v,i}$
%are irreducible cuspidal representations.

\medskip

We are now in the same situation as in \cite{boyer-aif} where we prove that
the conjecture 5.4.3 implies the conjecture 5.2.1 and
the translation in terms of automorphic congruences explained at the end of \S 5.2
The situation here is much more simple as $s=1$.

\begin{coro} \label{coro-congruence}
Let $\Pi$ be an irreducible automorphic representation of $G(\Am)$ which is $\xi$-comological of level $K$
and such that 
\begin{itemize}
\item its modulo $l$ Satake's parameters are given by $\mathfrak m$,

\item and its local component $\Pi_v$ at $v$ is isomorphic to 
$\Pi_v \simeq \st_{h}(\chi_{v}) \times \Psi_v,$
where $\chi_{v}$ is a characters and 
$\Psi_v$ is an irreducible representation of $GL_{d-h}(F_v)$.
\end{itemize}
Consider then any character $\chi'_{v}$ of $F_v^\times$ which is congruent to $\chi_{v}$ 
modulo $l$. Then there exists an irreducible automorphic representation $\Pi'$ of $G(\Am)$ 
which is $\xi$-cohomological of the same level $K$ and such that
\begin{itemize}
\item its modulo $l$ Satake's parameters are given by $\mathfrak m$,

\item its local component at $v$ is of the following shape
$$\Pi'_v \simeq \st_{h}(\chi'_{v}) \times \Psi'_v.$$
%where 
%\begin{itemize}
%\item for $i=2,\cdots,r$, the character $\chi'_{v,i}$ is congruent to $\chi_{v,i}$ modulo $l$
%and inertially equivalent to $\chi'_{v,1}$;
%
%\item $\Psi'_v \simeq \st_{s_1}(\pi'_{v,1}) \times \cdots \times \st_{s_u}(\pi'_{v,u})$ has
%a cuspidal support disjoint from the Zelevinsky line associated to $\chi'_{v,1}$.
%\end{itemize} 
\end{itemize}
\end{coro}
%
%Moreover if $\Psi_v \simeq  \st_{t_1}(\chi_{v,1}) \times \cdots \times \st_{t_r}(\chi_{v,r})$
%for some characters $\chi_{v,1},\cdots,\chi_{v,r}$, then $\Psi'_v$ is isomorphic to some 
%$\st_{t_1}(\chi'_{v,1}) \times \cdots \times \st_{t_r}(\chi'_{v,r})$
%where for $i=1,\cdots,r$, the character $\chi'_{v,i}$ is congruent to $\chi_{v,i}$ modulo $l$.
%
%\rem As explained before, using the main results of \cite{boyer-duke}, 
%one can replace in the previous
%proposition, the characters $\chi_{v}$ and $\chi'_{v}$ by irreducible cuspidal representations
%$\pi_{v}$ and $\pi'_{v}$ such that their common modulo $l$ reduction remains supercuspidal.
%
%\begin{proof}
%The existence of $\Pi'$ of level $K^p$ follows from (\ref{eq-equality}) and the description 
%of the cohomology in the previous proposition. For the remaining, i.e. the level at $p$
%and the case when $\Psi_v \simeq  \st_{t_1}(\chi_{v,1}) \times \cdots \times \st_{t_r}(\chi_{v,r})$,
%it follows trivially from the previous proof of the theorem  \ref{theo-principal}.
%
%
%\end{proof}

\bibliographystyle{plain}
\bibliography{bib-ok}

\end{document}